\documentclass[12pt]{article}
\usepackage{amsmath,amsfonts,amssymb,amscd}
\usepackage{graphics}
\usepackage{graphicx}
\usepackage{subfigure}
\usepackage[T2A]{fontenc}

\newtheorem{theo}{Theorem}
\newtheorem{lem}{Lemma}[section]
\newtheorem{prop}{Proposition}[section]
\newtheorem{cor}{Corollary}[section]
\newtheorem{rem}{Remark}[section]

\makeatletter \@addtoreset{equation}{section} \makeatother

\newcommand{\mC}{\mathbb{C}}

\newcommand{\mR}{\mathbb{R}}
\newcommand{\mT}{\mathbb{T}}
\newcommand{\mZ}{\mathbb{Z}}
\newcommand{\mN}{\mathbb{N}}

\newcommand{\bc}{{\bf c}}
\newcommand{\bh}{{\bf h}}
\newcommand{\bj}{{\bf j}}
\newcommand{\bu}{{\bf u}}

\newcommand{\AAA}{{\cal A}}
\newcommand{\CC}{{\cal C}}
\newcommand{\DD}{{\cal D}}
\newcommand{\FF}{{\cal F}}
\newcommand{\HH}{{\cal H}}
\newcommand{\II}{{\cal I}}

\newcommand{\MM}{{\cal M}}
\newcommand{\PP}{{\cal P}}

\newcommand{\eps}{\varepsilon}
\newcommand{\ph}{\varphi}

\newcommand{\rank}{\operatorname{rank}}

\newcommand\qed{{\unskip\nobreak\hfil\penalty50
  \hskip2em\hbox{}\nobreak\hfil\mbox{\rule{1ex}{1ex} \qquad}
    \parfillskip=0pt \finalhyphendemerits=0\par\medskip}}

\begin{document}

\title
{Lagrangian tori near resonances of near-integrable
Hamiltonian systems}
\author{Medvedev A. G.\footnote{Lomonosov Moscow State University}, Neishtadt A. I.\footnote{Space Research Institute of RAS, Department of Mathematical Sciences, Loughborough University}, Treschev D. V.\footnote{Steklov Mathematical Institute of RAS}}
\maketitle

\begin{abstract}
In this paper we study families of Lagrangian tori that appear in a neighborhood of a resonance of a near-integrable Hamiltonian system. Such families disappear in the ``integrable'' limit $\eps\to 0$. Dynamics on these tori is oscillatory in the direction of the resonance phases and rotating with respect to the other (non-resonant) phases.

We also show that, if multiplicity of the resonance equals one, generically these tori occupy a set of large relative
measure in the resonant domains in the sense that the relative
measure of the remaining ``chaotic'' set is of order $\sqrt\eps$.
Therefore for small $\eps > 0$ a random initial condition in a $\sqrt \eps$-neighborhood of a single resonance occurs inside this set (and therefore generates a quasi-periodic motion) with a probability much larger than in the ``chaotic'' set.

We present results of numerical simulations and discuss the form of projection of such tori to the action space.
\end{abstract}

\section{Projection of trajectories to the action space}

Consider a symplectic map
\begin{eqnarray}
\label{map_sym}
& (y,x) \mapsto (y_+,x_+), \qquad
  y\in\mR^2, \quad x\in\mT^2  &\\
\label{map_form}
& y_+ = y - \eps\, \partial V / \partial x, \quad
  x_+ = x + y_+, \qquad
  V = V(x) . &
\end{eqnarray}
For $\eps = 0$ the dynamics is integrable and $(y,x)$ are action-angle variables on the phase space $\mR^2\times\mT^2$. We choose for definiteness
$$
  V = a_1\cos(x_1 + \ph_1) + a_2\cos(x_2 + \ph_2)
       + a_3\cos(x_1 - x_2 + \ph_3)
$$
with constant $a_1,a_2,a_3,\ph_1,\ph_2,\ph_3$. We study trajectories numerically having in mind the idea to observe some interesting images and then find theoretic explanations. To visualize objects (closures of trajectories) lying in the 4-dimensional phase space, we project them to the action plane $\mR^2 = \{y\}$. We fix $a_j\sim 1$, $\ph_j$, $\eps\sim 1/10$, take a random initial condition and see.

\begin{figure}
\centering
\vskip-2mm
\subfigure
{\includegraphics[scale=0.75]{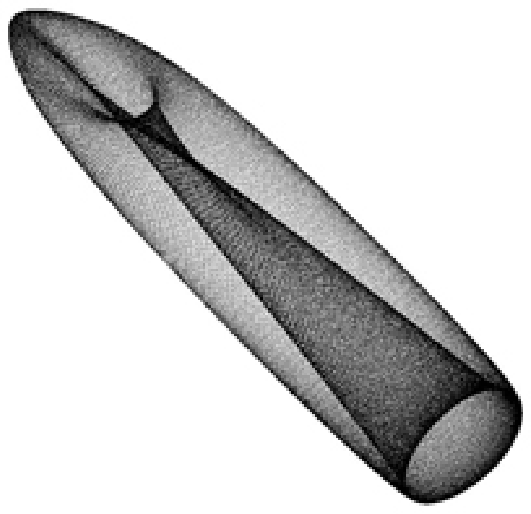}}
\subfigure
{\includegraphics[scale=0.75]{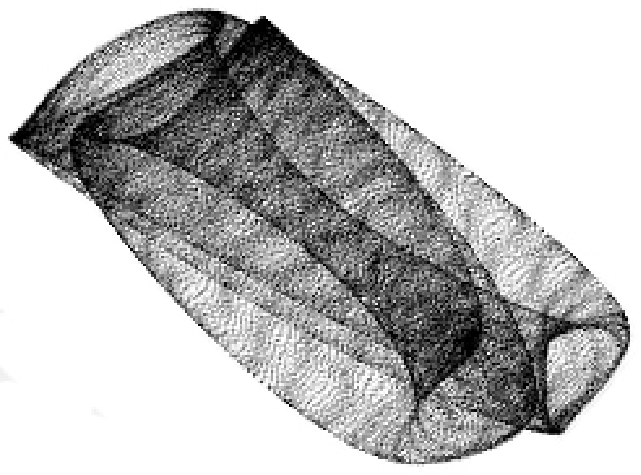}}
\subfigure
{\includegraphics[scale=0.75]{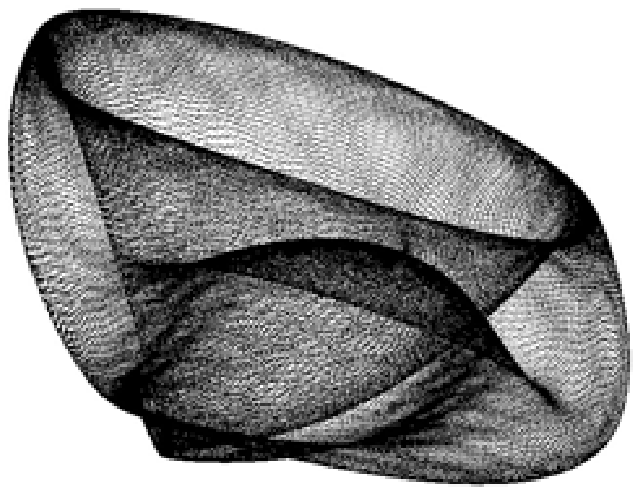}}
\caption{KAM-tori}
\label{fig:torus1}
\end{figure}

Typical images are presented on Fig. \ref{fig:torus1}. They can be easily identified with ordinary KAM-tori.  We see typical singularities (folds and pleats) of Lagrangian projections (these singularities are presented for example, in \cite{AKN}).

\begin{figure}
\centering
\vskip-2mm
\subfigure
{\includegraphics[scale=0.75]{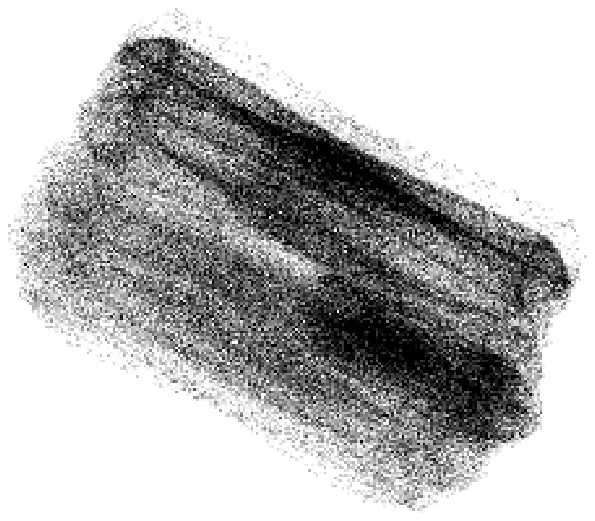}}
\subfigure
{\includegraphics[scale=0.75]{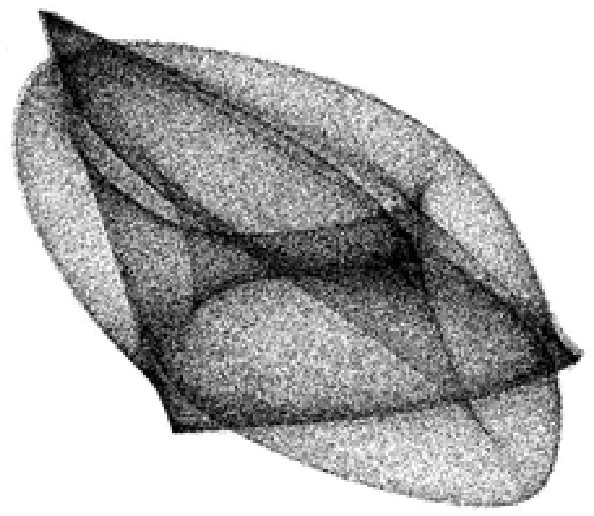}}
\subfigure
{\includegraphics[scale=0.25]{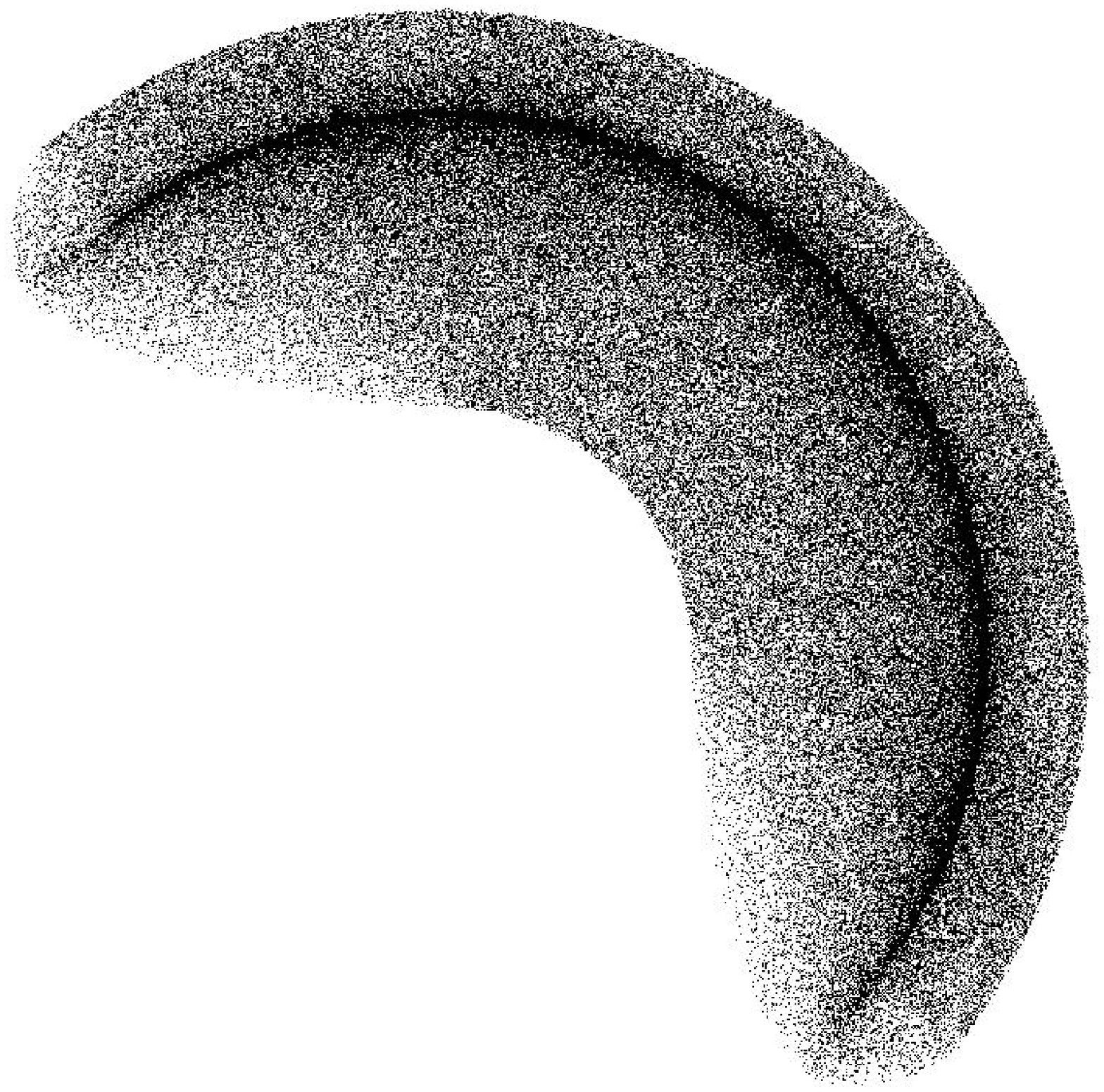}}
	\caption{Chaotic trajectories: near KAM-tori (center,right) and
             near a strong resonance (left)}
	\label{fig:chaotic}
\end{figure}

It is also easy to observe chaotic trajectories, Fig.~\ref{fig:chaotic}. The trajectories presented in Fig. \ref{fig:chaotic}, right and center, are situated near KAM-tori, but chaotic effects create a small ``defocusing'' in comparison with tori from Fig. \ref{fig:torus1}.

\begin{figure}
\centering
\vskip-2mm
\subfigure
{\includegraphics[scale=0.75]{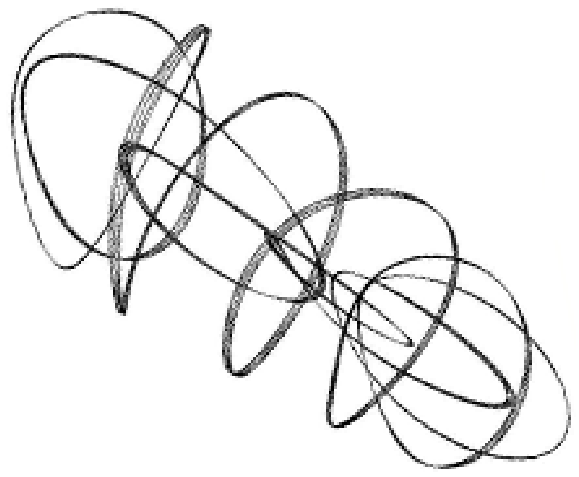}}
\subfigure
{\includegraphics[scale=0.75]{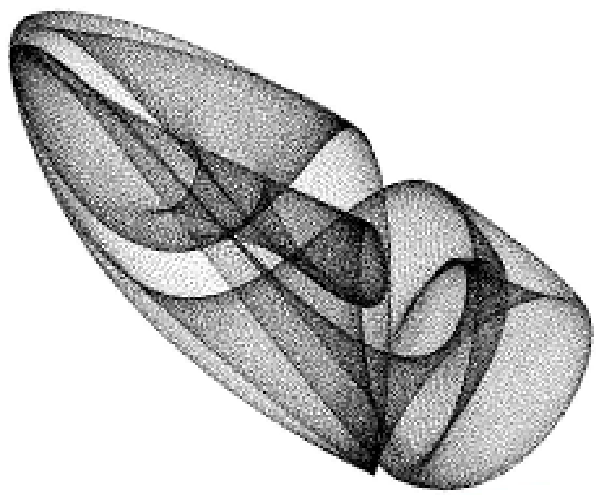}}
\subfigure
{\includegraphics[scale=0.75]{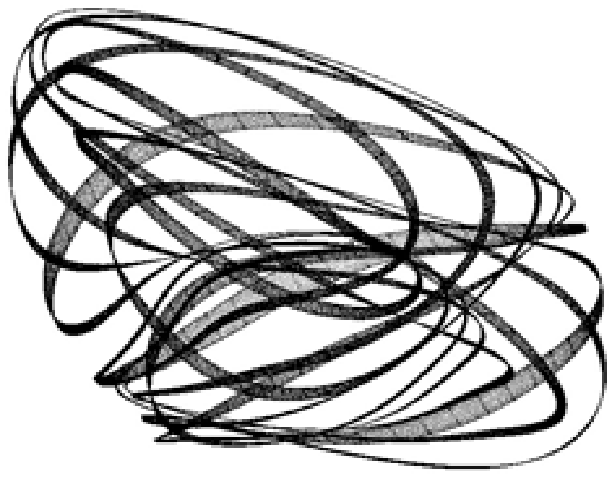}}
	\caption{Closed "ribbons"}
	\label{fig:ribbon}
\end{figure}

\begin{figure}
\centering
\vskip-2mm
\subfigure
{\includegraphics[scale=1]{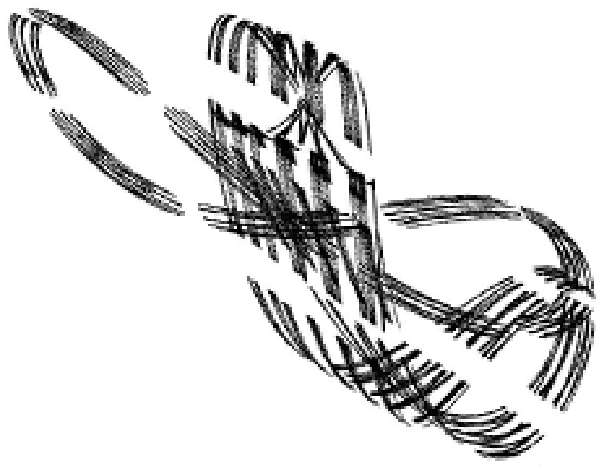}}
\subfigure
{\includegraphics[scale=1]{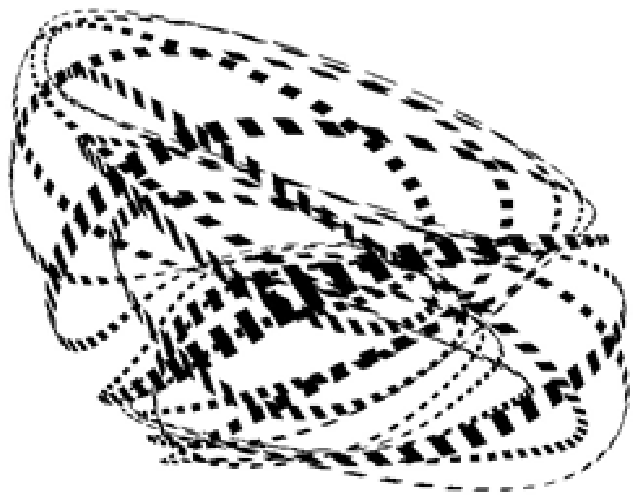}}
 \caption{Another type of quasi-periodic trajectories, ``sequences of spots''}
 \label{fig:spots}
\end{figure}

Trying more initial conditions, it is possible to obtain more mysterious objects, Fig. \ref{fig:ribbon}, which look like closed ribbons. Further attempts lead to more exotic images, Fig.~\ref{fig:spots}, looking as finite sequences of spots.

In this paper we discuss a mechanism which generates such structures. In particular we show that these objects form sets of positive measure in the phase space. Hence probability to observe them is also positive.

\section{Lower-dimensional tori near resonances}

As usual it is more convenient to write formulas for flows although numerics are faster, simpler and more precise for maps. Consider a real-analytic near-integrable Hamiltonian system
\begin{eqnarray}
\label{ham_sys}
&  \dot X = \partial H / \partial Y, \quad
   \dot Y = -\partial H / \partial X, \qquad
   Y\in\DD\subset\mR^N, \quad
   X\in\mT^N. \\
\label{ham_iniini}
&  H = H_0(Y) + \eps H_1(Y,X) + O(\eps^2), \qquad
   \eps\ge 0. &
\end{eqnarray}
Below we use the following notation for such a system:
$$
  (P, \omega_P, H), \qquad
  P = \DD\times\mT^N, \quad
  \omega_P = dY\wedge dX,
$$
where the symplectic manifold $(P,\omega_P)$ is the phase space.

The map (\ref{map_sym})--(\ref{map_form}) can be regarded as the Poincar\'e map for some system (\ref{ham_sys})--(\ref{ham_iniini}) with $N=3$ on an energy level $M_h = \{H=h=\mbox{const}\}$ (see for example \cite{PT}). Hence the dimension of the phase space $6=2N$ drops by 1 because of the reduction to $M_h$ and by another 1 because of the passage to a (hyper)surface $\Sigma\subset M_h$ transversal to the flow. Below $N\ge 3$ is arbitrary.

The vector $\nu(Y) = \partial H_0 / \partial Y$ is called an unperturbed frequency. For a fixed $Y=Y^0$ we have a fixed frequency $\nu^0 = \nu(Y^0)$. Any equation
\begin{equation}
\label{resonan}
  \langle\nu^0,K\rangle = 0, \qquad K\in\mZ^N\setminus\{0\}
\end{equation}
is called a resonance. The word ``resonance'' is also attributed to the integer vector $K$, satisfying (\ref{resonan}).

Given a constant $\nu^0\in\mR^N$ all the corresponding resonances $K$ (together with $0\in\mZ^N$) form a resonance $\mZ$-module
$$
     g(\nu^0)
   = \{K\in\mZ^N : K \mbox{ is a resonance}\} \cup \{0\}.
$$
If $g(\nu^0)=0$, the frequency vector $\nu^0$ is said to be nonresonant. We define
$$
  l = l(\nu^0) = \rank(g(\nu^0)), \quad
  m = m(\nu^0) = N-l,
$$
where $\rank(g(\nu^0))$ is the number of generators in $g(\nu^0)$. Informally speaking, $l$ is the number of independent resonances for the frequency vector $\nu^0$.
Invariant torus is called resonant (nonresonant) if the dynamics on the torus is quasi-periodic with a resonant (nonresonant) frequency vector.
Any torus $$\mT^N_{Y^0} = \{(Y,X) : Y=Y^0\}$$ is invariant with respect to the unperturbed flow
$$
  (Y,X) \mapsto (Y,X+\nu(Y)t), \qquad  t\in\mR.
$$
Then $\mT^N_{Y^0}$ is foliated by the $m$-tori
$$
    \mT^m_{Y^0,X^0}
  = \mbox{closure}
      \big(\big\{ (Y,X) : Y=Y^0,\; X=X^0+\nu^0 t,\; t\in\mR
      \big\}\big).
$$
Computing frequency vector corresponding to the unperturbed quasi-periodic motion on $\mT^m_{Y^0,X^0}$, it is easy to show that the tori $\mT^m_{Y^0,X^0}$ are non-resonant.

Note that in the non-resonant case $(l=0)$ we have
$\mT^m_{Y^0,X^0} = \mT^m_{Y^0}$ for any $X^0\in\mT^N$. If $l>0$ then $m\ne N$ and the foliation is non-trivial.

If $l>0$ then a generic perturbation destroys $\mT^m_{Y^0}$, \cite{Poi}. However generically some tori $\mT^m_{Y^0,X^0}$ survive a perturbation even if $\nu^0$ is resonant. To present the corresponding result, we fix a $\mZ$-module $g\subset\mZ^n$. Consider the resonance set
$$
  \Sigma_g = \{Y\in\DD : g(\nu(Y)) = g\}.
$$
Under natural non-degeneracy conditions\footnote
{the functions $\langle K^{(j)},\partial/\partial Y\rangle H_0(Y)$ are independent in $\DD$ where $K^{(1)},\ldots,K^{(l)}$ are generators of $g$} $\Sigma_g$ is a real-analytic submanifold,
$\dim\Sigma_g = m$.

It is convenient to study system (\ref{ham_sys})--(\ref{ham_iniini}) in a $\sqrt\eps$-neighborhood of the torus $\mT^N_{Y^0}$, $Y^0\in\Sigma_g$ by using the scaling
$$
  Y = Y^0 + \sqrt\eps\widetilde Y, \quad
  X = \widetilde X, \qquad
      H(Y,X,\eps)
    = H_0(Y^0) + \sqrt\eps\widetilde H(\widetilde Y,\widetilde X,\sqrt\eps).
$$
Then the system $(P,\omega_P,H)$ turns to the system
$(\widetilde P,\omega_{\widetilde P},\widetilde H)$,
\begin{eqnarray*}
&   \widetilde H
  = \langle \nu^0,\widetilde Y \rangle
   + \frac{1}{2}\sqrt\eps
       \langle \Pi\widetilde Y,\widetilde Y \rangle
   + \sqrt\eps\, H_1(Y^0,\widetilde X)
   + O(\eps), & \\
&  \widetilde P = \mR^N\times\mT^N, \quad
   \omega_{\widetilde P} = d\widetilde Y\wedge d\widetilde X, \quad
   \Pi = H_{0\,YY}(Y^0). &
\end{eqnarray*}

The tori $\mT^m_{Y^0,X^0}$ which survive the perturbation are generated by fixed points of some Hamiltonian system which is obtained from the initial one by using averaging, neglecting some higher order perturbative terms, and reduction of the order. Now we turn to description of these steps.

For any function $f = \sum_{K\in\mZ^N} f^K e^{i\langle K,X\rangle}$
consider the averaging
\begin{equation}
\label{averaging}
    \langle f\rangle_g
  = \sum_{K\in g} f^K e^{i\langle K,X\rangle}.
\end{equation}
Hence, $\langle\cdot\rangle_g$ is a projector, removing all nonresonant Fourier harmonics.

Now it is natural to perform the following standard coordinate change:
$$
  (\widetilde Y,\widetilde X) \mapsto (\widehat Y,\widehat X),
  \qquad
    \widetilde Y
  = \widehat Y + \sqrt\eps\,\widetilde S_{\widetilde X}, \quad
    \widetilde X
  = \widehat X,
$$
where $\widetilde S = \widetilde S(\widetilde X)$ is a solution of the (co)homological equation
\begin{equation}
\label{cohomol}
     \langle \nu^0,\widetilde S_{\widetilde X}\rangle
   + H_1(Y^0,\widetilde X)
 =   v(Y^0,\widetilde X), \qquad
     v(Y,X)
 =   \langle H_1(Y,X) \rangle_g .
\end{equation}
Under standard Diophantine conditions a real-analytic solution of equation (\ref{cohomol}) exists and unique up to a $g$-invariant additive term $\langle \widetilde S \rangle_g$. For example, $\widetilde S$ can be chosen so that
$\langle\widetilde S\rangle_g = 0$.

In the new coordinates we have the system
$(\widehat P,\omega_{\widehat P},\widehat H)$,
\begin{eqnarray*}
&\displaystyle
 \widehat P = \mR^N\times\mT^N, \quad
 \omega_{\widehat P} = d\widehat Y\wedge d\widehat X, \quad
 \widehat H = H^g + O(\eps), & \\
&\displaystyle
    H^g
 =  \langle\nu^0,\widehat Y\rangle
  + \frac{\sqrt\eps}{2} \langle\Pi\widehat Y,\widehat Y\rangle
  + \sqrt\eps\, v(Y^0,\widehat X). &
\end{eqnarray*}

Consider the approximate system
$(\widehat P,\omega_{\widehat P},H^g)$, which is usually called the partially averaged system. Any critical point $X^0$ of the ``potential'' $v(Y^0,X)$ generates an invariant torus $\mT_{Y^0,X^0}^m$. To study linearization of the averaged system on $\mT_{Y^0,X^0}^m$ it is convenient to consider the reduced system. First we recall general invariant construction (see \cite{AKN}) and then give a more explicit coordinate form.

The system $(\widehat P, \omega_{\widehat P}, H^g)$ admits the symmetry group $G\cong\mT^m$. The Lie algebra associated with $G$ is naturally identified with
$$
    g_\perp
  = \{Q\in\mR^N : \langle Q,K\rangle = 0
                    \mbox{ for any $K\in g$}\}.
$$
Note that $\nu^0\in g_\perp$ and $\rank g_\perp = m$.

Action of $G$ on $P$ is Poissonian and the corresponding momentum map $\MM_G : P\to g_\perp^*$ is as follows:
$$
  \MM_G(\widehat Y,\widehat X)\, Q = \langle Q,\widehat Y \rangle
  \quad
  \mbox{for any } Q\in g_\perp .
$$
Reduction with respect to $G$ means
\smallskip

(a) fixing values of the first integrals
    $\langle Q,\widehat Y \rangle = \langle Q,\gamma\rangle$ for any $Q\in g_\perp$, where $\gamma = \mbox{const} \in\mR^N$,

(b) passage to the quotient phase space
$$
    \PP
  = \widehat P_\gamma / G, \qquad
    \widehat P_\gamma
  = \big\{(\widehat Y,\widehat X)\in\widehat P :
             \langle Q,\widehat Y\rangle
             = \langle Q,\gamma\rangle \;
                         \mbox{ for all } Q\in g_\perp \big\}.
$$
\smallskip

The reduced phase space has a canonical symplectic structure $\omega_{\gamma}$ \cite{AKN} while the Hamiltonian $\HH^g : \PP\to\mR$ is determined by the commutative diagram
$$
  \begin{array}{rcl}
  \widehat \PP_{\gamma} &  \stackrel{\mbox{\rm\small pr}}{\longrightarrow}  &  \PP \\
  {}_{H^g|_{\widehat \PP_{\gamma}}}\searrow\!\!\!\! &&\!\!\!\! \swarrow {}_{\HH^g}\\
      &  \mR          &
  \end{array}
$$
where $\mbox{pr}$ is the natural projection. The torus $\mT^m_{Y^0,X^0}$ turns to a fixed point $p$ in the reduced system $(\PP,\omega_{\gamma},\HH^g)$. Flow of the system $(P,\omega_P,H^g)$ near $\mT^m_{Y^0,X^0}$ is essentially determined by linear approximation of $(\PP,\omega_{\gamma},\HH^g)$ at $p$. To obtain this approximation, we turn to the coordinate form of the above order reduction.

Let $\Theta$ be an $N\times l$ matrix formed by the integer vectors $K^{(1)},\ldots,K^{(l)}$, generators of $g$. This matrix is not unique: for any $M\in SL(l,\mZ)$ (an integer $l\times l$ matrix with unit determinant) one may take $\Theta M$ instead of $\Theta$. We assume that the quadratic form determined by $\Pi$ is non-degenerate on $g_{\mR}$, where $g_{\mR}\subset\mR^N$ is the natural extension of $g(\nu^0)$. In other words, $\Theta^T\Pi\Theta$ is a non-degenerate $l\times l$-matrix.

Then we can take as local coordinates on $\PP$ the variables $\eta\in\mR^l$, $\xi\in\mT^l$ such that
$$
  \widehat Y = \Theta(\eta-\eta^0) + \gamma, \quad
  \xi = \Theta^T X, \qquad
  \eta^0 = (\Theta^T\Pi\Theta)^{-1} \Theta^T\Pi\gamma.
$$
Here the constant $\eta^0$ is chosen for convenience to remove from $\HH^g$ a term linear in $\eta$.

The form $\omega_{\gamma}$, the fixed point $p$, and the function $\HH^g$ are as follows:
$$
  \omega_{\gamma} = d\eta\wedge d\xi, \quad
  \eta(p) = 0, \; \xi(p) = \Theta^T X^0, \quad
    \HH^g
  =  \frac{\sqrt\eps}{2} \langle\Theta^T\Pi\Theta\eta,\eta\rangle
    + \sqrt\eps\, v_{\PP,Y^0}(\xi)
    + h_\gamma.
$$
Here $v_{\PP,Y^0} : \mT^l\to\mR$ is the unique function, satisfying the identity\footnote
{Explicit formula for $v_{\PP}$ is as follows. Let the Fourier expansion for $v$ be
$$
    v(Y,X)
  = \sum_{K\in g} v^K(Y) e^{i\langle K,X\rangle}
  = \sum_{j\in\mZ^l} v^{\Theta j}(Y)
                           e^{i\langle\Theta j,X\rangle}.
$$
Then $v_{\PP,Y}(\xi) = \sum_{j\in\mZ^l} v^{\Theta j}(Y)
                               e^{i\langle j,\xi\rangle}$.}
$$
  v_{\PP,Y^0}(\Theta^T X) = v(Y^0,X).
$$
The constant
$h_\gamma = \frac12\langle\Pi(\gamma-\Theta\eta^0),
                             (\gamma-\Theta\eta^0)\rangle$
can be ignored.

\begin{theo} \emph{(\cite{LiYi}).}
\label{theo:poi++}
Suppose that $\det\Pi\ne 0$ in $D$. Then for any sufficiently small $\eps>0$ there exists a set $\Lambda_\eps\subset\Sigma_g$ such that for each $Y^0\in\Lambda_\eps$ and for each nondegenerate critical point $\xi^0$ of $v_{\PP,Y^0}$ the perturbed system admits a real-analytic invariant $m$-torus $\mT^m_{Y^0,\xi^0}(\eps)$. This torus is close to $\mT^m_{Y^0,X^0}$, where $X^0$ is any point satisfying the equation $\xi^0 = \Theta^T X^0$. Moreover, $\mT^m_{Y^0,\xi^0}(\eps)$ carries a quasi-periodic motion with the same frequency vector.

The perturbed invariant $m$-tori constitute a finite number of $m$-parameter Whitney smooth families. The relative Lebesgue measure of $\Lambda_\eps$ on the surface
$$
  \{ Y\in\Sigma_g : v_{\PP,Y^0}
     \mbox{ has nondegenerate critical points}\}
$$
tends to 1 as $\eps\to 0$.
\end{theo}

Hamiltonian of the linear approximation for
$(\PP,\omega_{\gamma},\HH^g)$ at the fixed point $p$ is
$$
    \HH^g_{lin}
  = \frac{\sqrt\eps}{2} \langle\Theta^T\Pi\Theta\eta,\eta\rangle
   + \frac{\sqrt\eps}{2} \langle V\zeta,\zeta\rangle, \qquad
    V
  = \frac{\partial^2 v_{\PP}}{\partial\xi^2} (0), \quad
    \zeta
  = \xi - \xi^0 .
$$
Therefore eigenvalues $\pm\lambda_1,\ldots,\pm\lambda_l$ of the fixed point $\eta=\zeta=0$ in the system $(\PP,\omega_{\gamma},\HH^g_{lin})$ satisfy the equation
$$
  \det (\eps\Theta^T\Pi\Theta V + \lambda^2) = 0.
$$
If all $\lambda_j$ are purely imaginary, $p$ and the corresponding torus $\mT^m_{Y^0,\xi^0}(\eps)$ are said to be normally elliptic. The ``opposite'' case is normally hyperbolic, where no $\lambda_j$ is purely imaginary.

For $m=1$ Theorem \ref{theo:poi++} was proven by Poincar\'e \cite{Poi}. In this case no small divisors appear and the proof is  based on the ordinary implicit function theorem. The equation $m=N$ corresponds to the (ordinary) KAM-theorem for Lagrangian tori.

Hyperbolic case with arbitrary $m$ is presented in \cite{Tre_hyp}.
In \cite{CKLY} the case of arbitrary $m$ and arbitrary normal behavior of the perturbed tori is treated. In \cite{LiYi} it is shown that an additional condition (the so called $g$-nondegeneracy of $H_0$), introduced in \cite{CKLY}, can be skipped. A statement  analogous to Theorem \ref{theo:poi++} should be true in infinite dimension but as far as we know this has not been proven yet.

Since for non-trivial $g$ all the tori $\mT^m_{Y^0,\xi^0}(\eps)$ are lower-dimensional ($m<N$), their total measure in the phase space $D\times\mT^N$ vanishes. In other words they are practically invisible in numerical experiments. This does not mean that they are inessential for dynamics. For example, hyperbolic lower-dimensional tori and their asymptotic manifolds are known as elementary links which form transition chains, forming a basis for the Arnold diffusion.

\section{Visible objects}
\label{sec:visible}

In this paper we are interested in visible objects. More precisely, in invariant tori of dimensional $N$. Geometry of their projections to the action space depends on the order of a resonance at which these tori appear.
\medskip

{\bf No resonance}. For example, such objects are ordinary ($N$-dimensional) KAM-tori. If $\eps>0$ is small, KAM tori form a large Cantorian set: the measure of the complement $\CC(\eps)$ to this set in $D\times\mT^N$ does not exceed a quantity of order $\sqrt\eps$, \cite{AKN,Laz,Nei,Poe,SW}. The measure estimates of $\CC(\eps)$ for degenerated systems are contained in \cite{Arn,Mazz,Vaid,RandP}.

In the first approximation in $\sqrt\eps$ projection of a KAM-torus to the action space $\DD$ has the form $\widehat Y = 0$ under the condition $l = \rank g = 0$. In the original coordinates we have:
$$
  \big\{ Y\in\DD :
         Y = Y^0
             + \sqrt\eps\widetilde S_{\widetilde X}(\widetilde X),
         \quad \widetilde X\in\mT^N
  \big\},
$$
where $\widetilde S$ satisfies (\ref{cohomol}).

The discrete system (\ref{map_sym})--(\ref{map_form}) can be obtained from a Hamiltonian system (\ref{ham_sys})--(\ref{ham_iniini}) with $N=3$ degrees of freedom on an energy level $H=\mbox{const}$ by passing to the Poincar\'e map on the section $\{X_3=0\}$. Hence, the first in $\sqrt\eps$ approximation of the objects presented in Fig. \ref{fig:torus1} are sets of the form
\begin{equation}
\label{proj_tor}
  \big\{ (Y_1,Y_2) :
         Y_j = Y_j^0
          + \sqrt\eps\widetilde S_{\widetilde X_j}
                   (\widetilde X_1,\widetilde X_2,0),\quad
         j = 1,2, \quad
         (\widetilde X_1,\widetilde X_2)\in\mT^2
  \big\}.
\end{equation}
For a random initial condition the probability to occur on one of such torus is greater than $1- C_0 \sqrt \eps$ for some $C_0 > 0$.

The ``resonance'' set $\CC(\eps)$ contains other families of quasi-periodic motions. Here it is reasonable to distinguish the case of a single resonance ($\rank g = 1$) and the case of a multiple resonance ($\rank g > 1$).
\medskip

{\bf Single resonance}. In this case dimension of the commutative symmetry group $G$ is $N-1$. Therefore the system $(P,\omega_P,H^g)$ is completely integrable. Informally speaking, it is a product of $N-1$ ``rotators'' and a ``pendulum'' (such representation in a small neighborhood of a single resonance is discussed in \cite{ARN2}), where the pendulum is determined by the reduced system with one degree of freedom $(\PP,\omega_{\gamma},\HH^g)$. The variables $\eta$ and $\xi$ are 1-dimensional while $\Theta = (\Theta_1,\ldots,\Theta_N)^T$ is just an integer vector.

Solutions on which the pendulum motions are rotations are ordinary KAM-tori while solutions on which the pendulum oscillates, lie in $\CC(\eps)$. When we say that the pendulum oscillates, we mean that in the system $(\PP,\omega_{\gamma},\HH^g)$ the angular variable $\xi$ changes periodically in an interval $J\subset\mT^1$, $J\ne\mT^1$. In the first approximation in $\sqrt\eps$ the corresponding trajectory of (\ref{map_sym})--(\ref{map_form}) fills the set
$$
  \big\{ (Y_1,Y_2) :
         Y_j = Y_j^0
          + \sqrt\eps\widetilde S_{\widetilde X_j}
                   (\widetilde X_1,\widetilde X_2,0),\quad
         j = 1,2, \quad
        \Theta_1\widetilde X_1 + \Theta_2\widetilde X_2\in J,\quad
        (\widetilde X_1,\widetilde X_2)\in\mT^2
  \big\}.
$$
This is a ``ribbon-like'' subset of (\ref{proj_tor}). This explains the structure of sets in Fig. \ref{fig:ribbon}.

Since the system $(P,\omega_P,H^g)$ is integrable, the set of phase points in a $\sqrt\eps$-neighborhood of the torus $\mT^N_{Y^0}$ lying outside invariant $N$-tori of the original system $(P,\omega_P,H)$, has a small relative measure. Precise statement, Theorem \ref{theo:ccc}, is given in Section \ref{sec:single}. Hence if an initial condition is taken randomly the probability to obtain a quasi-periodic orbit like Fig. \ref{fig:ribbon} is not less than $C_1 \sqrt \eps, \quad C_1 > 0$, because the width of a resonance domain corresponding to a single resonance is $\sim \sqrt \eps$. Note that pictures analogous to Fig. \ref{fig:ribbon} can be found in \cite{Simo}.
\medskip

{\bf Multiple resonance}. If $\rank g > 1$, the systems $(P,\omega_P,H^g)$ and $(\PP,\omega_{\gamma},\HH^g)$ are generically non-integrable. Therefore existence of invariant $l$-tori in the latter one is not straightforward provided the energy $\HH^g$ is not very big or not very small. A standard source for such tori is a neighborhood of a totally elliptic fixed point. However, a totally elliptic fixed point may not exist if the ``kinetic energy''
$\frac12 \sqrt\eps\langle\Theta^T\Pi\Theta\eta,\eta\rangle$ is indefinite: a simple example is
$$
  \HH^g = \sqrt\eps\, (\eta_1\eta_2 + \cos\xi_1 + \cos\xi_2).
$$

System (\ref{map_sym})--(\ref{map_form}) corresponds to a positive definite kinetic energy and trajectories presented in Fig. \ref{fig:spots} present nonlinear versions of small oscillations near a totally elliptic periodic orbit in the corresponding system with 3 degrees of freedom. A random initial condition lies on one of such tori with probability of order  $\eps$, because the measure of a resonant domain corresponding to a double resonance is of order $\eps$. Unlike the case of a single resonance only a small portion of this domain is filled with tori, in general.

\section{Invariant $N$-tori at a single resonance}
\label{sec:single}

Putting $N=n+1$, consider the Hamiltonian system (\ref{ham_sys})--(\ref{ham_iniini}) in a neighborhood of the resonance $\Sigma_g\times\mT^{n+1}$, where $g$ is generated by the vector $K^0\in\mZ^{n+1}\setminus \{0\}$ with relatively prime components. Hence, we plan to study invariant tori, located in the vicinity of a single resonance
\begin{equation}
\label{resonance}
    \Sigma
  = \Sigma_g
  = \{Y\in \DD : \langle K^0,\nu(Y)\rangle = 0\}.
\end{equation}
We assume that the unperturbed system is non-degenerate and $K^0$ is not a light-like vector:
\begin{eqnarray}
\label{non-deg}
& \mbox{det}\, H''_{0YY} \ne 0, \qquad
  Y\in \DD, & \\
\label{res_nondeg}
&        \AAA(Y)
  \equiv \big\langle K^0, H''_{0\, YY} \, K^0
         \big\rangle
  \ne 0. &
\end{eqnarray}
Note that (\ref{res_nondeg}) means that for any $Y$ satisfying (\ref{resonance}) the function
$\lambda\mapsto H_0(Y - \lambda K^0)$ has a non-degenerate
critical point $\lambda = 0$.

Now our aim is to give a defenition of the oscillatory part of the resonance domain and to introduce convenient notation for the main KAM theorem.

By (\ref{res_nondeg}) the resonant set $\Sigma\subset\DD$ is a smooth hypersurface transversal to the constant vector field $K^0$. The equation
\begin{equation}
\label{ddlam}
         \frac{d}{d\lambda} H_0(Y - \lambda K^0)
  \equiv \big\langle K^0,\nu(Y - \lambda K^0) \big\rangle
    =    0
\end{equation}
has a real-analytic solution $\lambda = \lambda(Y)$ in $\DD$ near
$\Sigma$.

We have a smooth map $\chi : U(\Sigma)\to\Sigma$, where
$U(\Sigma)$ is a neighborhood of the resonance $\Sigma$ and
$\chi(Y) = Y - \lambda(Y) K^0$.

Let $\langle\cdot\rangle^{K^0}\equiv\langle\cdot\rangle_g$ be the operator of resonant averaging
$$
  f = \sum_{K\in\mZ^{n+1}} f^K\, e^{i\langle K,X\rangle}
  \;\mapsto\;
    \langle f\rangle^{K^0}
  = \sum_{j=-\infty}^\infty
      f^{j K^0}\, e^{ijq}, \qquad
    q = \langle K^0,X \rangle .
$$

Consider the Hamiltonian system $(P,\omega_P,H_{K^0})$,
\begin{equation}
\label{ham_ave}
    H_{K^0}(Y,X)
  = H_0(Y) + \eps\,\bu(Y,q), \quad
    \bu(Y,q)
  = \big\langle H_1\big(\chi(Y), X\big)
    \big\rangle^{K^0} .
\end{equation}
The function $\bu$ is $2\pi$-periodic in the resonant phase $q$ and
$$
  \bu(Y + \lambda K^0,q) = \bu(Y,q) \quad
  \mbox{for any $\lambda$ in a neighborhood of $0\in\mR$}.
$$
In fact, below the quantity $\lambda(Y) K^0 = Y - \chi(Y)$ (some sort of distance to the resonant surface) will be of order $\sqrt\eps$.

Since
$   H_0(Y)
 =  H_0(\chi(Y))
   + \frac12 \AAA(\chi(Y))\, \lambda^2(Y)
   + O(\lambda^3(Y)),
$
Hamiltonian (\ref{ham_ave}) can be presented in the form
\begin{equation}
\label{HH_K0}
    H_{K^0}(Y,X)
 =  H_0(\chi(Y))
   + \frac12 \AAA(\chi(Y))\, \lambda^2(Y)
   + \eps\bu(Y,q)
   + O(\lambda^3(Y)).
\end{equation}

For any vector $J\in\mR^{n+1}$ such that $\langle K^0,J\rangle = 0$ the function $\langle Y,J\rangle$ is a first integral. Therefore the system $(P,\omega_P,H_{K^0})$ is completely integrable. It is responsible for the dynamics of the original system
near the resonance (\ref{resonance}). Below we only deal with the
oscillatory part ${\cal D}_{os}$ of the resonance domain, where
${\cal D}_{os}\subset\DD\times\mT^{n+1}$ is defined as follows. Let
$q_{\min}(Y)$ and $q_{\max}(Y)$ be  points of global minimum and maximum of $\bu(Y,q)$ for fixed $Y\in\Sigma$:
\begin{equation}
\label{minmax}
  \bu(Y,q_{\min}(Y)) = \min_{q\in\mT} \bu(Y,q), \quad
  \bu(Y,q_{\max}(Y)) = \max_{q\in\mT} \bu(Y,q).
\end{equation}
Then we define
\begin{eqnarray}
\nonumber
 \!\!\!\!\!\!\!
     \DD_{os}\!\!
 &=&\!\! \big\{ (Y,X)\in\DD\times\mT^{n+1} \; : \\
\label{DD}
 && \quad
              \eps \bu(Y,q_{\min}(Y))
            < \frac12 \AAA(\chi(Y))\, \lambda^2(Y)
   + \eps\bu(Y,q)
            < \eps \bu(Y,q_{\max}(Y))
     \big\} .
\end{eqnarray}

If $\bu(Y,q)$ is not a constant as a function of $q$, the domain
$\DD_{os}$ belongs to  an $O(\sqrt\eps)$-neighborhood of the resonance
$\Sigma\times\mT^{n+1}$. On almost any orbit of this flow located in $\DD_{os}$ the resonant phase $q$ oscillates between two quantities $q_1$ and $q_2$, depending on initial conditions and such that
$$
  |q_2 - q_1| < 2\pi, \quad
  \dot q |_{q = q_1} = \dot q |_{q = q_2} = 0.
$$
These orbits lie on $(n+1)$-dimensional Lagrangian tori. Below we prove that for small values of $\eps$ all these tori except a set of a small measure survive the perturbation.

To formulate the result, for any $Y$ in a neighborhood of $\Sigma$ consider the Hamiltonian system
\begin{equation}
\label{hatbh}
  \dot q = \sqrt\eps\, \hat\bh'_p(\chi(Y),p,q), \quad
  \dot p = - \sqrt\eps\, \hat\bh'_q(\chi(Y),p,q).
\end{equation}
with one degree of freedom, the Hamiltonian
$$
    \hat\bh(\chi(Y),p,q)
  = \frac12 \AAA\big( \chi(Y) \big) p^2
       + \bu(Y,q)
$$
and the symplectic structure $\frac{1}{\sqrt\eps} dp\wedge dq$. The point $Y$ is regarded as a parameter. Recall that by (\ref{res_nondeg}) $\AAA\ne 0$. This system coincides with $(\PP,\omega_{\gamma},\HH^g)$ from Section \ref{sec:visible} written in slightly other terms.

\begin{figure}
\begin{center}
\includegraphics[scale = 1]{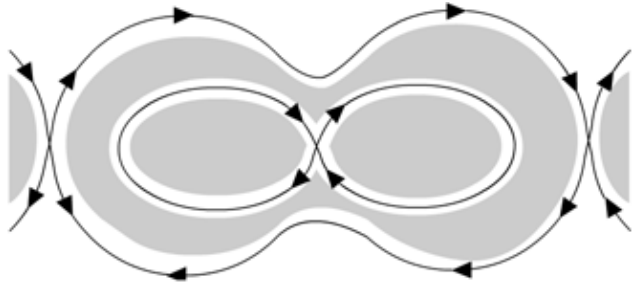}
\caption{Phase portrait of system (\ref{hatbh}).
         Section of the oscillatory domain $\DD_0$ by
         $\{\chi(Y) = \mbox{const}\}$ is marked grey.}
\label{fig:pp}
\end{center}
\end{figure}

\begin{prop}
\label{prop:ave}
The point $\chi(Y)$ is a constant of motion in the system $(P,\omega_P,H_{K^0})$.

The variables $\lambda(Y) = \sqrt\eps p$ and $q$ in the averaged system up to $O(\eps + \lambda^2)$ satisfy equations (\ref{hatbh}).
\end{prop}

Indeed, the first statement of Proposition \ref{prop:ave} follows from the relation $\dot Y\parallel K_0$.

Applying the operator
$\langle K^0,\partial/\partial Y\rangle$ to (\ref{ddlam}) we get:
$\langle K^0,\lambda_Y(Y) \rangle = 1+O(\lambda)$. Then by (\ref{HH_K0}) in system
$(P,\omega_P,H_{K^0})$
\begin{eqnarray*}
     \dot\lambda
 &=& \lambda'_Y\,\dot Y
  =  - \eps\bu'_q(Y,q) +O(\eps\lambda)
  = -\eps\hat\bh'_q + O(\eps\lambda), \\
     \dot q
 &=& \langle K^0, H'_{K^0\, Y} \rangle
  =  \AAA(\chi(Y))\,\lambda + O(\lambda^2)+O(\eps)
  =  \sqrt\eps\,\hat\bh'_p + O(\eps + \lambda^2).
\end{eqnarray*}
\qed

The projection
$\pi : U(\Sigma\times\mT^{n+1})\to\Sigma\times\mR\times\mT$,
$$
  \pi(Y,X) = (\chi(Y),p,q), \qquad
  p = \eps^{-1/2} \lambda(Y), \quad
  q = \langle K^0,X\rangle
$$
maps the domain $\DD_{os}$ to $\hat\DD_{os} = \pi(\DD_{os})$, see Fig. \ref{fig:pp}.

Let the closed curve $\gamma(Y,p,q)$ be the connected component of the set
$$
  \big\{(\tilde p,\tilde q) :
         \hat\bh(\chi(Y),\tilde p,\tilde q) = \hat\bh(\chi(Y),p,q)
  \big\}
$$
containing the point $(p,q)$. We define the action variable $I$ and the Hamiltonian function $\bh$:
\begin{equation}
\label{I=I}
  I = I(Y,p,q)
    = \frac{1}{2\pi}\int_{\gamma(Y,p,q)} \tilde p\,d\tilde q, \qquad
  \bh(\chi(Y),I) = \hat\bh(\chi(Y),p,q).
\end{equation}
Note that if the energy levels $\bh(\chi(Y),p,q) = \mbox{const}$ are not connected (i.e., consist of several curves $\gamma$), the function $\bh$ is not single-valued.

If $\gamma = \gamma(Y,p,q)$ is a closed smooth curve, the torus
\begin{equation}
\label{tor}
    \mT_{Y^0,\gamma}^{n+1}(\eps)
  = \big\{ (Y,X) : \chi(Y) = Y^0,
                   (\eps^{-1/2}\lambda(Y),q)\in\gamma
    \big\}
\end{equation}
is invariant for the system $(P,\omega_P,H_{K^0})$ up to terms of order $O(\eps+\lambda^2)$. For small $\eps>0$ majority of tori (\ref{tor}) survive the perturbation and exist in the original system. Any surviving torus has to satisfy several additional conditions.
\smallskip

(1) The frequency vector associated with $\mT_{Y^0,\gamma}^{n+1}$ is Diophantine.
\smallskip

This is a standard assumption which holds on almost all tori.
\smallskip

(2) The system $(P,\omega_P,H_{K^0})$ is nondegenerate on
$\mT_{Y^0,\gamma}^{n+1}$.
\smallskip

This condition essentially means that $\bh$ is real-analytic in a neighborhood of $\mT_{Y^0,\gamma}^{n+1}$ and
\begin{equation}
\label{twist}
   |\bh'_{I}| < c, \quad |\bh''_{II}| < c,\quad |\bh'_{I}|^{-1} < c, \quad |\bh'_{II}|^{-1} < c.
\end{equation}
Therefore we have to replace $\DD_{os}$ by a smaller domain $\DD_0$ by throwing out a small neighborhood of asymptotic manifolds (where the tori degenerate) and a small neighborhood of tori on which the twist conditions (\ref{twist}) is violated. Let $\mu$ be the measure in the phase space $\DD\times\mT^{n+1}$ generated by the symplectic structure $\omega_P$. Then the measures $\mu(\DD)$ and $\mu(\DD_0)$ are both of order $O(\sqrt\eps)$.

\begin{theo}
\label{theo:ccc}
Suppose that the system
$(P,\omega_P,H)$ is real-analytic. If $\eps_0>0$ is sufficiently small, then for all positive $\eps \le \eps_0$ any Diophantine torus $\mT_{Y^0,\gamma}^{n+1}(\eps)\subset\DD_0$ survives the perturbation.
For some constant $C>0$ independent of $\eps$, measure of union of such tori in $\DD_0$ is not less than $\mu(\DD_0) - C \eps$.
\end{theo}

\begin{rem}
Note, that $\mu(\DD) - \mu(\DD_0) \sim \sqrt \eps$. To obtain $\mu(\DD) - \mu(\DD_0) \sim \eps$  we need to replace in a small neighborhood of asymptotic manifolds the twist conditions (\ref{twist}) by a weaker ones. The strong proof is not ready now.
\end{rem}

\section{Preliminaries}
Beginning from this place up to end of the  paper we prove Theorem \ref{theo:ccc}.
\begin{itemize}
\item{All vectors by default are regarded as columns.
      For any $u\in\mR^m$ and any $m\times m$-matrix $A$
      we use the notation
$$
    |u|
  = \max_{1\le j\le m} |u_j|, \quad
    |A|
  = \max_{0\ne u\in\mR^m} \frac{|Au|}{|u|} .
$$
 The brackets $\langle\,,\rangle$ denote the standard Euclidean
 scalar product: $\langle u,v \rangle = \sum_{j=1}^m u_j v_j$.
}

\item{$\mu$ denotes the standard Lebesgue measure on $\mR^m$.}

\item{Prime denotes a partial derivative e.g.,
      $f'_{y_k} = \partial f / \partial y_k$. If $y\in\mR^m$, and $f:\mR^m\to\mR$ then
      $f^\prime_y$ is regarded as a vector and
      $f^{\prime\prime}_{yy}$ as a matrix.}

\item{For any function $f:\mT^m\to\mR$ we define its average
 \begin{equation}
 \label{average}
    \langle f \rangle
  = \frac{1}{(2\pi)^m} \int_{\mT^m} f(x)\, dx.
 \end{equation}
 The same notation is used if $f$ depends on other variables. In
 this case to avoid misunderstanding we use for (\ref{average}) the
 notation $\langle\;\rangle_x$.
 }

\item{Below $c_1,c_2,\ldots$ denote positive constants.
      If $c_j$ depends on another constant, say, $\alpha$,
      we write $c_j(\alpha)$. Dependence on the dimension
      $n$ is not indicated.}
\end{itemize}

To present the system $(P,\omega_P,H^g)$ in a form convenient for application of KAM procedure, we have to perform several preliminary coordinate changes.

{\bf (a)}. Consider a matrix $M\in GL(n+1,\mZ)$ such that $K^0$ is its last column. In the new coordinates
$$
  \hat Y = M^{-1} Y, \quad
  \hat X = M^T X
$$
resonance (\ref{resonance}) takes the form
\begin{equation}
\label{resonance1}
         \hat\nu_{n+1}(\hat Y)
 \equiv  H'_{0\,\hat Y_{n+1}}(\hat Y)
    =   0, \qquad
        H_0(\hat Y)
    =   H_0(Y).
\end{equation}

To have more convenient coordinates in a neighborhood of this
resonance, we solve the first equation (\ref{resonance1}) with respect to $\hat Y_{n+1}$. This can be done locally because by
(\ref{res_nondeg})
$H''_{0\,\hat Y_{n+1}\hat Y_{n+1}} \ne 0$. We denote the result
$$
  \hat Y_{n+1} = G(\hat Y_1,\ldots,\hat Y_n),\qquad
  \hat\nu_{n+1}\big( \hat Y_1,\ldots,\hat Y_n,
                        G(\hat Y_1,\ldots,\hat Y_n)
               \big) \equiv 0.
$$

{\bf (b)}. Consider the change of the variables
\begin{eqnarray*}
& y = (y_1,\ldots,y_n) = (\hat Y_1,\ldots,\hat Y_n), \quad
  p = \eps^{-1/2} (\hat Y_{n+1} - G(y)), & \\
& x = (x_1,\ldots,x_n)
    = \big( \hat X_1 + G'_{y_1} q,\ldots,\hat X_n + G'_{y_n} q
      \big), \quad
  q = \hat X_{n+1}.
&
\end{eqnarray*}
We are interested in motions which are oscillatory in the coordinate $q$. Therefore it is not necessary to assume periodicity of this change with respect to $q$. Below $q$ lies in an interval while the variables $x$ are still angular: $x\in\mT^n$.

The resonance $\Sigma$ in the new coordinates locally takes the form $\{p=0\}$. In particular, $y$ are local coordinates on $\Sigma$.

Then the symplectic structure and the Hamiltonian take the form
\begin{eqnarray}
\nonumber
 &   \omega
  =  dy\wedge dx + \sqrt\eps\, dp\wedge dq, & \\
\label{H=Lam}
 &   \Lambda(y)
    + \eps \big( \frac12 A(y) p^2 + u(y,q) \big)
    + \eps U_1(y,x,q)
    + \eps^{3/2} U_2(y,p,x,q,\sqrt\eps) , &
\end{eqnarray}
where
$$
    \Lambda(y)
  = H_0(y,G(y)), \quad
    A(y)
  = H''_{0\,\hat Y_{n+1}\hat Y_{n+1}}(y,G(y)),
$$
$u,U_1,U_2$ are real-analytic, and average of $U_1$ with respect to $x$ vanishes: $\langle U_1 \rangle_x = 0$.

By (\ref{res_nondeg}),(\ref{ham_ave})
\begin{equation}
\label{u=u}
      A(y)
  =  \AAA\big( \chi(Y) \big) ,\quad
      u(y,q)
  =   \bu(\chi(Y),q).
\end{equation}
Neglecting the terms $\eps U_1 + \eps^{3/2} U_2$, we obtain an
integrable system which can be regarded as a skew-product of an
$n$-dimensional rotator in variables $y,x$ and a (generalized)
pendulum in variables $p,q$.

We put (see (\ref{minmax}))
$$
  \hat q_{\min}(y) = q_{\min}(\chi(Y)), \quad
  \hat q_{\max}(y) = q_{\max}(\chi(Y)).
$$
Then the domain $\DD_{os}$ (see (\ref{DD})) takes the form
\begin{equation}
\label{DDnew}
     \DD_{os}
  =  \big\{ (y,p,x,q) : u(y,\hat q_{\min})
                       < \frac12 A(y) p^2 + u(y,q)
                       <  u(y,\hat q_{\max})
     \big\}.
\end{equation}

{\bf (c)}. Let $W(y,I,q)$ be a generating function which introduces
action-angle variables $I,\ph$ in domain (\ref{DDnew}) for the system with one degree of freedom, the symplectic structure
$dp\wedge dq$ and the Hamiltonian
$\frac12 A(y) p^2 + u(y,q)$ for any values of the parameters $y$:
$$
    p = W'_q, \quad
  \ph = W'_{I}, \quad
  h(y,I) = \frac12 A(y) p^2 + u(y,q).
$$
Then the canonical change with the generating function
$\hat y x + \sqrt\eps W(\hat y,I,q)$
$$
       y = \hat y, \quad
  \hat x = x + \sqrt\eps W'_{\hat y},\quad
    p = W'_q, \quad
  \ph = W'_{I}
$$
transforms the symplectic structure to
$d\hat y\wedge d\hat x + \sqrt\eps dI\wedge d\ph$ and Hamiltonian
(\ref{H=Lam}) to
\begin{equation}
\label{ham_ini}
     \Lambda(\hat y)
    + \eps h(\hat y,I)
    + \eps \hat V(\hat y,I,\hat x,\ph,\sqrt\eps) , \quad
     \langle \hat V\rangle_x = O(\sqrt\eps),
\end{equation}
where the functions $h,\hat v,\hat V = U_1 + \sqrt \eps U_2$ are real-analytic. The
function $h$ satisfies the equation
$$
  h(y,I) = \bh(\chi(Y),I), \qquad
  Y\in U(\Sigma),
$$
where $\bh$ is defined in (\ref{I=I}).
\medskip
Below we skip hats for brevity.

\section{Initial KAM Hamiltonian}

For any set $D \subset\mR^{n+1}$ let $U_a(D)\subset\mC^{n+1}$ be the following neighborhood:
$$
     U_a(D)
  =  \{ (y+\eta,I+ \zeta): (y,I) \in D,
                         |\eta| < \sqrt \eps a, |\zeta| < a
     \}, \qquad
    \eta\in\mC^n, \; \zeta\in\mC.
$$
For any function $f : D\to\mR$ which admits a real-analytic extension to $U_a(D)$ we put
$$
    |f|_{a}
  = \sup_{z \in U_a(D)} |f(z)|.
$$
Note that this norm is anisotropic in $y$ and $I$ directions.

Let $U_{b}(\mT^{n+1})$ be a complex neighborhood of $\mT^{n+1}$
$$
    U_{b}(\mT^{n+1})
  = \{ (x + \xi,\ph + \kappa):
      (x,\ph)\in\mT^{n+1},\, |\xi| < b, |\kappa| < b\},\qquad
    \xi\in\mC^n, \; \kappa\in\mC.
$$

For any function $f:\mT^{n+1}\to\mR$ which admits a real-analytic extension to $U_b(\mT^{n+1})$ we put
$$
    |f|_{b}
  = \sup_{z \in U_b(\mT^{n+1})} |f(z)|.
$$

For functions, real analytic on $D\times\mT^{n+1}$ we define
$|f|_{a,b}$ as the corresponding double supremums over
$$
    U_{a,b}(D\times\mT^{n+1})
  = U_a(D)\times U_b(\mT^{n+1}).
$$

Consider the Hamiltonian system with the symplectic structure
\begin{equation}
\label{sym_str}
  dy\wedge dx + \sqrt\eps\, dI\wedge d\ph
\end{equation}
and the real-analytic Hamiltonian (see (\ref{ham_ini}))
$$
    H_0
  =  \Lambda(y)
    + \eps h_0(y,I)
    + \eps v_0(y,I,\ph,\sqrt\eps)
    + \eps u_0(y,I,x,\ph,\sqrt\eps) , \qquad
    v_0(y,I,\ph,0) = 0, \quad
    \langle u_0\rangle_x = 0,
$$
where $v_0 = \langle V \rangle_x$, $u_0 = V - \langle V \rangle_x$ and the points $(y,I,x,\ph)$ lie in a complex neighborhood
$$
   U_{a_0,b_0}(D_0\times\mT^{n+1}), \qquad
  D_0 \subset \mR^{n+1}
$$
for some $a_0,b_0 > 0$.

The above analyticity assumptions mean that there exist
$\eps_0,\bar s_0,s_0 > 0$ such that for any
$0 \le \eps < \eps_0$
\begin{equation}
\label{0}
      |\Lambda|_{a_0}
  \le \bar s_0, \quad
      |h_0|_{a_0}
  \le c_h, \quad
      |u_0|_{a_0,b_0}
  \le s_0, \quad
      |v_0|_{a_0,b_0}
  \le \sqrt\eps\bc.
\end{equation}
Assumptions of Theorem \ref{theo:ccc} imply the following non-degeneracy conditions:
\begin{eqnarray}
\label{Lam-1}
&\!\!\!\!\!\!
    \underline{c}_\Lambda \le  |\det \Lambda''_{yy}| \le \overline{c}_\Lambda, \quad |\Lambda''_{yy}|_{a_0}
  \le c_\Lambda , \quad
      |\Lambda^{\prime\prime\,-1}_{yy}|_{a_0}
  \le c_\Lambda , \quad
      |h'_{0\, I}|_{a_0}
  \le c'_h, \quad
  		|h'_{0\, I}|^{-1}_{a_0}
  		 \le c'_h, \quad
   & \\
\label{h0''}
&\!\!\!\!\!\!
      |h''_{0\, II}|_{a_0}
  \le c'_h , \quad
        |h''_{0\, II}|^{-1}_{a_0}
  \le c'_h , \quad
      |\sqrt \eps h''_{0\, Iy}|_{a_0}
  \le c''_h , \quad
      |\eps h^{\prime\prime}_{0\, yy}|_{a_0}
  \le c''_h .  &
\end{eqnarray}
Let $\eps$ be sufficiently small, then we can assume that $c''_h$ is small because $c''_h \sim \sqrt \eps$. Below we assume that $c''_h \le c''_{h\,0}(\underline{c}_{\Lambda},c_{\Lambda},c'_h)$.

\section{The Hamiltonian $H_m$}

Below all functions depend smoothly on $\eps$. For
brevity we do not write $\eps$ in their arguments.

As usual KAM procedure includes a converging sequence $\FF_0,\FF_1,\ldots$ of coordinate changes and a converging sequence of Hamiltonians
$H_0,H_1,\ldots$
$$
  \FF_m : U_{a_{m+1},b_{m+1}}(D_{m+1}\times\mT^{n+1})
        \to U_{a_m,b_m}(D_m\times\mT^{n+1}), \quad
  H_{m+1} = H_m \circ \FF_m.
$$

Consider an increasing sequence $\{N_j\in\mZ\}$ ($N_j$ is the maximal order of a resonance essential on the $j$-th step), a
decreasing sequence $\{\lambda_j > 0\}$
($\sqrt \eps \,\lambda_j$ determines the width of resonance strips on the $j$-th step) and the function $\bj : \mN\to\mN$ defined by the inequality (let $N_{-1} = 0$)
\begin{equation}
\label{N<N}
  N_{\bj(r) - 1} < r \le N_{\bj(r)} \quad
  \mbox{for all $r > 0$}.
\end{equation}
Then $\bj(r)$ is the number of the first step on which the resonance of order $r$ is essential.

Consider two positive decreasing sequences $a_m,b_m$,
\begin{equation}
\label{ab}
  a_m = a_{m+1} + 6\sigma_m, \quad
  b_m = b_{m+1} + 3\delta_m
\end{equation}

Suppose that on the $m$-th step we have the Hamiltonian
\begin{equation}
\label{Hm}
    H_m
  = \Lambda(y)
    + \eps h_m(y,I)
    + \eps v_m(y,I,\ph)
    + \eps u_m(y,I,x,\ph), \qquad \langle u_m\rangle_x = 0.
\end{equation}

The function $H_m$ is defined in a complex neighborhood $U_{a_m,b_m}(D_m\times\mT^{n+1})$
\begin{equation}
\label{D=D}
      D_{m+1}
   =  D_m\setminus
      \bigcup_{|K|\le N_m,\,k\ne 0}
       U_{a_m}(Q_{K,m}) , \qquad
      K
  =  \Big(\begin{array}{c} k \\ k_0
         \end{array}\Big)
 \in \mZ^{n+1} , \quad
     k\in\mZ^n,
\end{equation}
where the resonant strips $Q_{K,m}$ are defined with the help of the sequences $N_j$ and $\lambda_j$:
\begin{eqnarray}
\label{Q}
&     Q_{K,m}
   =  \big\{ (y,I) \in U_{a_m}(D_m) :
             \big| \langle\nu_m(y,I),K\rangle
             \big|
         \le \lambda_{\bj(|K|)} (1+2^{-m-1}) \sqrt \eps
      \big\}, & \\
\label{psi}
&    \nu_m(y,I)
  =  \Big( \begin{array}{c}
             \Lambda'_y(y) + \eps h'_{m\, y}(y,I) \\
             \sqrt\eps\, h'_{m\,I}(y,I)
           \end{array}
     \Big) . &
\end{eqnarray}
\begin{prop}
\label{prop:B}
For any
$m = 1,2,\ldots$
$$
      \mu(D_0\setminus D_m)
  \le c_\mu \sqrt \eps,
$$
where $c_{\mu} > 0$ is independent of $\eps$.
\end{prop}

{\bf Inductive assumptions}. For $(y,I)\in U_{a_m}(D_m)$ the
following estimates hold:
\begin{eqnarray}
\label{||m<}
 \!\!\!\!\!\!\!\!\!&    |v_m|_{a_m,b_m}
  \le s_m, \quad
      |u_m|_{a_m,b_m}
  \le s_m , & \\
\label{h_I}
 \!\!\!\!\!\!\!\!\!&    |h_m|_{a_m}
  \le (2 - 2^{-m}) c_h ,  \quad
     |h'_{m\,I}|_{a_m}^{-1}
  \le (2 - 2^{-m}) c'_h ,  \quad
       |h''_{m\, II}|_{a_m}^{-1}
  \le (2 - 2^{-m}) c'_h , & \\
\label{h_y}
\!\!\!\!\!\!\!\!\!&\displaystyle
      |h''_{m\,II}|_{a_m}
  \le (2 - 2^{-m}) c'_{h} , \quad
      |\sqrt \eps h''_{m\,Iy}|_{a_m}
  \le (2 - 2^{-m}) c''_{h} , \quad
      |\eps h''_{m\,yy}|_{a_m}
  \le (2 - 2^{-m}) c''_{h} . &
\end{eqnarray}

\section{The KAM-step}

For any natural $N$ and a periodic function
$$
  f:\mT^{n+1}\to\mR, \qquad
     f(x,\ph)
  = \sum_{K = (k,k_0)\in\mZ^{n+1}}
        f^K e^{i\langle k,x\rangle + ik_0\ph}
$$
we define the cut off
\begin{equation}
\label{cutoff}
    \Pi_N f(x,\ph)
  = \sum_{|K|\le N,\, k\ne 0}
         f^K e^{i\langle k,x\rangle + ik_0\ph}.
\end{equation}
Then by Lemma \ref{lem:cutoff} for any real-analytic $f$ such that $|f|_b < \infty$ and for
any $\delta\in (0,b)$
\begin{equation}
\label{PIN}
      \big| f - \langle f \rangle_x - \Pi_N f
      \big|_{b-\delta}
  \le \frac{C}{\delta} \Big(N + \frac{1}{\delta}\Big)^n
            e^{-N\delta}\, |f|_{b} .
\end{equation}
\smallskip

By using Hamiltonian (\ref{Hm}), we introduce the canonical\footnote
{i.e., preserving symplectic structure (\ref{sym_str})}
change of variables
$(y,I,x,\ph) \mapsto (\hat y,\hat I,\hat x,\hat\ph)$,
determined by the generating function
$\hat yx
  + \sqrt\eps\big(\hat I\ph
                  + S(\hat y,\hat I,x,\ph)\big)$:
$$
       y = \hat y + \sqrt \eps S^\prime_x, \quad
       I = \hat I +  S^\prime_\ph, \quad
  \hat x = x + \sqrt \eps S^\prime_{\hat y}, \quad
 \hat\ph = \ph + S^\prime_{\hat I} ,
$$
where the arguments $(\hat y,\hat I)$ are supposed to lie in $ U_{a_m - \sigma_m}(D_{m+1})$.
\begin{rem}
\label{rem:inD_m+1}
$|\cdot|_{a_m-\sigma,b_m -\delta}$ (resp. $|\cdot|_{a_m-\sigma}$) denotes the norm in $U_{a_m-\sigma,b_m -\delta}(D_{m+1}\times\mT^{n+1})$ (resp. $U_{a_m-\sigma}(D_{m+1})$).
\end{rem}

By definition the function $S$ is a solution of the homological equation
\begin{equation}
\label{homologic}
    \Big\langle \nu_m(\hat y,\hat I),
                \Big(\!\!\begin{array}{c}
                       S'_x \\ S'_\ph
                     \end{array}\!\!
                \Big)(\hat y,\hat I,x,\ph)
    \Big\rangle
  = - \sqrt \eps \,\Pi_{N_{m}} V_m(\hat y,\hat I,x,\ph),
    \qquad
    V_m = v_m + u_m .
\end{equation}

\begin{prop}
\label{prop:homolo}
For any $m = 0,1,\ldots$ there exists a solution of (\ref{homologic}) where
\begin{equation}
\label{L}
      |S|_{a_m,b_m - \delta_m}
  \le L_m s_m,\qquad
      L_m
   =  \sum_{j=0}^{m} 2\frac{N_j^{n+1}}{\lambda_j}
         e^{- (n+1) \delta_m N_{j-1}}.
\end{equation}
\end{prop}

The Hamiltonian (\ref{Hm}) takes the form
\begin{equation}
\label{Hm+1}
    \widetilde H_m
  =  \Lambda(\hat y)
    + \eps h_m(\hat y,\hat I)
    + \eps v_m(\hat y,\hat I,\hat\ph)
    + \eps\tilde v_m(\hat y,\hat I,\hat\ph)
    + \eps\tilde u_m(\hat y,\hat I,\hat x,\hat\ph) , \qquad
      \langle\tilde u_m \rangle_x = 0 .
\end{equation}

\begin{rem}
\label{rem:ybary}
Below we show that that
\begin{eqnarray}
\label{ang}
&     |S'_x|_{a_m - \sigma_m,
                                     b_m - 2\delta_m}
  \le \sigma_m, \quad
      |S'_\ph|_{a_m - \sigma_m,b_m - 2\delta_m}
  \le \sigma_m, & \\
\label{act}
&     |\sqrt \eps S'_y|_{a_m - \sigma_m,
                                    b_m - 2\delta_m}
  \le \delta_m, \quad
      |S'_I|_{a_m - \sigma_m,b_m - 2\delta_m}
  \le \delta_m. &
\end{eqnarray}
Estimates (\ref{ang})--(\ref{act}) imply that the coordinate change is well-defined for
$(\hat y,\hat I,\hat x, \hat \ph)
 \in U_{a_m - \sigma_m,b_m - 2\delta_m}
    (D_{m+1}\times \mT^{n+1})$.
\end{rem}

\begin{prop}
\label{prop:RRRR} For $m\ge 1$ estimates (\ref{ang})--(\ref{act}) imply the inequalities
\begin{eqnarray}
\label{|v|}
 &    |\tilde v_m + \tilde u_m|_{a_m - \sigma_m, b_m - 2\delta_m}
  \le  \tilde s_m,  & \\
\label{tilde s_m}
 &  \displaystyle{
      \tilde s_m
   =  (c_{\Lambda} + c'_h)\Big(\frac{L_m s_m}{\delta_m}\Big)^2
     + (n+2) \frac{L_m s_m^2}{\sigma_m\delta_m} + \frac{C}{\delta_m} \Big(N_{m} + \frac{1}{\delta_m}\Big)^n
       e^{-N_{m}\delta_m}\, s_m }.
   \end{eqnarray}
\end{prop}

\section{An additional step}

Consider the symplectic transformation
$(\hat I,\hat\ph)\mapsto (\bar I,\bar\ph)$ with generating function
$\bar I\hat\ph + \sqrt\eps \widetilde S(\bar y,\bar I,\hat\ph)$ which introduces action-angle variables in the system with one degree of freedom and Hamiltonian
\begin{equation}
\label{hm+1}
    h_m(\hat y,\hat I) + v_m(\hat y,\hat I,\hat\ph)
  = h_{m+1}(\bar y,\bar I)  .
\end{equation}
The variables $\hat y = \bar y$ are regarded as parameters. We
extend this map to a canonical transformation of the whole phase
space:
\begin{equation}
\label{^to-}
  \hat y = \bar y, \quad
  \hat I = \bar I + \widetilde S'_{\hat\ph}, \quad
  \bar x = \hat x + \sqrt\eps \widetilde S'_{\bar y},\quad
  \bar\ph = \hat\ph + \widetilde S'_{\bar I}.
\end{equation}
Then Hamiltonian (\ref{Hm+1}) takes the form
$$
     H_{m+1}
  =  \Lambda(\bar y)
    + \eps h_{m+1}(\bar y,\bar I)
    + \eps v_{m+1}(\bar y,\bar I,\bar\ph)
    + \eps u_{m+1}(\bar y,\bar I,\bar x,\bar\ph) , \qquad
    \langle u_{m+1}\rangle_{\bar x} = 0.
$$

\begin{prop}
\label{prop:Hn+1}
Suppose that
\begin{equation}
\label{cs}
      |v_0|_{a_0} \le \sqrt\eps \bc
  \le \frac{\sigma_0 \delta_0}{2 c'_h}, \quad
      |v_m|_{a_m-\sigma_m,b_m-2\delta_m}
  \le  s_m, \quad
       s_m    \le \frac{\sigma_m \delta_m}{4c'_h}, \qquad
  m \ge 1.
\end{equation}
Then for any $m\ge 0$
\begin{eqnarray}
\label{tildeS}
|\tilde S'_{\hat\ph}|_{a_m - 4\sigma_m,b_m - 2\delta_m} \le 8c'_h s'_m , \quad
 |\tilde S|_{a_m - 4\sigma_m,b_m - 2\delta_m} \le 16\pi c'_h s'_m ,\\
\label{|h-h|}
|h_{m+1} - h_m|_{a_m - 4\sigma_m,b_m - 2\delta_m} \le \frac{8 c_h c'_h s'_m}{\sigma_m} , \quad
  |u_{m+1} + v_{m+1}|_{a_{m+1},b_{m+1}}\le  s'_{m},
 \end{eqnarray}
where $s'_0 = \sqrt\eps \bc$ and $s'_m = s_m$ for all $m\ge 1$.
\end{prop}

\begin{rem}
\label{rem:addstep}
Below we show that
\begin{eqnarray}
\label{AddStepang}
&    |\widetilde S'_{\bar I}|_{a_{m+1},b_{m+1}}
 \le \delta_m, \quad
     |\sqrt\eps\widetilde S'_{\bar y}|_{a_{m+1},
                                        b_{m+1}}
  \le \delta_m, & \\
\label{AddStepact}
&   |\widetilde S'_{\hat\ph}|_{a_{m+1},b_{m+1}}
 \le \sigma_m.
\end{eqnarray}
Estimates (\ref{AddStepang}),(\ref{AddStepact}) imply that the coordinate change (\ref{^to-}) is well-defined for
$(\bar y,\bar I,\bar x, \bar \ph) \in
 U_{a_{m+1},b_{m+1}}(D_{m+1}\times \mT^{n+1})$.
\end{rem}

\section{The sequences
         $\sigma_m,\delta_m,s_m,L_m,N_m,\lambda_m$}

We define $\sigma_m$ and $\delta_m$ by (\ref{ab}) and put
\begin{eqnarray}
\label{defsigmadelta}
    &&\sigma_m
  = \frac{a'_0}{6} \, 2^{-(2n+5)(m+1)}, \quad \delta_m
  = \frac{b_0}{3} \, 2^{-(m+1)}, \quad  a'_0 = \frac{a_0}{2^{2n+5} -1}\\
\label{defsNLambda}
        &&s_m
 \; =\; s_0\, e^{- c_s m - 2^m},\quad
     N_m
 = c_N \, \, 2^{2m}, \quad
        \lambda_m
 \; =\; c_\lambda\, 2^{-(2n+3)m}.
\end{eqnarray}
The constants $a_0$, $b_0$, $s_0$ are a priori fixed. We can choose only $c_N$, $c_s$, $c_{\lambda}$ and $\eps$. First we fix $c_s$, then we define $c_N = c_N(c_s)$ and $c_{\lambda} = c_{\lambda}(c_N,c_s)$. Below we describe the process of this choosing.

\begin{prop}
\label{prop:L}
Suppose that the sequences $\delta_m,N_m$, and $\lambda_m$ are
defined by (\ref{defsigmadelta}) and (\ref{defsNLambda}). Then the
sequence $L_m$, defined by (\ref{L}) satisfies the estimate
\begin{equation}
\label{L<}
      L_m
  \le \frac{c_N^{n+1}}{c_{\lambda}}2^{(4n+5)(m+2)}.
\end{equation}
\end{prop}

To show that our choice of the sequences $a_m,b_m,\sigma_m,\delta_m,s_m,L_m,N_m,\lambda_m$ makes the procedure converging, we have to check that our assumptions
(\ref{||m<})--(\ref{h_y}), (\ref{ang}), (\ref{act}) and (\ref{AddStepang}), (\ref{AddStepact}) hold. The remaining part of this section contains this check.
\medskip

\subsection{Several estimates}
By using (\ref{defsigmadelta})--(\ref{L<}) we obtain:
\begin{eqnarray*}
\Big( \frac{L_m s_m}{ \delta_m} \Big)^2 &\le& 9 \frac{4^{8n+11} c_N^{2(n+1)}2^{(8n+12)m}s_0}{b_0^2 c_{\lambda}^2} s_0 e^{-2c_sm - 2^{m+1}} \le c_{Ls}2^{(8n+12)m} e^{-c_s(m-1)} s_{m+1},\\
  \frac{L_m s^2_m}{ \delta_m \sigma_m} &\le& 18 \frac{2^{10n+16} c_N^{n+1} 2^{(6n+11)m}s_0}{b_0 a'_0 c_{\lambda}}s_0 e^{-2c_sm - 2^{m+1}} \le c_{Ls}2^{(8n+12)m} e^{-c_s(m-1)} s_{m+1},
\end{eqnarray*}
where $c_{Ls} = \max\Big(9 \frac{4^{8n+11} c_N^{2(n+1)}s_0}{b_0^2 c_{\lambda}^2}, 18 \frac{2^{10n+16} c_N^{n+1}s_0 }{b_0 a'_0 c_{\lambda}}\Big)$.
For the third term of (\ref{tilde s_m}) we have
\begin{eqnarray*}
 &&\frac{C}{\delta_m} \Big(N_{m} + \frac{1}{\delta_m}\Big)^n e^{-N_{m}\delta_m}\, s_m \le \frac{6C}{b_0}\Big(c_N + \frac{6}{b_0}\Big)^n 2^{(2n+1)m} e^{-b_0 c_N 2^{m}/6}s_0 e^{-c_sm - 2^m}\\
	&&\le \frac{6C}{b_0}\Big(c_N + \frac{6}{b_0}\Big)^n e^{(2n + 1)m + c_s - (b_0c_N/6-2)2^{m} - 2^m}s_{m+1}.
\end{eqnarray*}

\subsection{Inequalities (\protect\ref{||m<})}
Rewrite (\ref{tilde s_m}) in terms of $s_{m+1}$
\begin{eqnarray*}
      \tilde s_m
   &\le&  (c_{\Lambda}+c'_h)c_{Ls}2^{(8n+12)m} e^{-c_s(m-1)} s_{m+1} + (n+2) c_{Ls}2^{(8n+12)m} e^{-c_s(m-1)} s_{m+1}\\
     && + \frac{6C}{b_0}\Big(c_N + \frac{6}{b_0}\Big)^n e^{(2n + 1)m + c_s - (b_0c_N/6-2)2^{m} - 2^m}s_{m+1} .
\end{eqnarray*}
If $c_s > c_{s_0} = 16n + 24$ and $c_{\lambda} > c_{\lambda\,0}(a_0,b_0,s_0,c_s,c_N,c_{\Lambda})$, then for all $m \ge 0$ we have
$$
c_{Ls}2^{(8n+12)m} e^{-c_s(m-1)} \le \frac{1}{3}\max((c_{\Lambda}+c'_h)^{-1},(n+2)^{-1}).
$$
For sufficently large $c_{N\,0}$ for all $c_N > c_{N\,0}(b_0,c_s)$ we obtain
$$
\frac{6C}{b_0}\Big(c_N + \frac{6}{b_0}\Big)^n e^{(2n + 1)m + c_s - (b_0c_N/6-2)2^{m} - 2^m} \le \frac{1}{3}.
$$
It implies $\tilde s_m \le s_{m+1}$. From Proposition \ref{prop:Hn+1} and the last inequality follows (\ref{||m<}).

\subsection{Inequalities (\protect\ref{h_I}),(\protect\ref{h_y})}

By (\ref{|h-h|}) $h_1 - h_0 = O(\sqrt\eps)$. Hence it is sufficient to check (\ref{h_I}) and (\protect\ref{h_y}) only for $m \ge 1$. For large $c_s > c_{s\,1} = c_{s\,1}(c'_h,a_0,s_0)$:
\begin{eqnarray*}
       |h_{m+1} - h_m|_{a_{m+1}+2\sigma_m}
 &\le& \frac{8 c_h c'_h s_m}{\sigma_m} \le 48 c_h \frac{2^{2n+5} c'_h s_0}{a'_0}e^{-(c_s - 2n-5) m - 2^m} \le c_h 2^{-m-1},\\
        |h''_{m+1\,II} - h''_{m\,II}|_{a_{m+1}}
 &\le& \frac{8 c_h c'_h s_m}{\sigma_m^3} \le 1728 \frac{2^{6n+15} c_h c'_h s_0}{a'^2_0} e^{-(c_s -6n-15)m - 2^m} \\
 &\le& c'_h 2^{-m-1}.
\end{eqnarray*}

Note, that for $c_s > c_{s\,2} = c_{s\,2}(c_h,c'_h,a_0,s_0)$ we have
\begin{eqnarray*}
&&|h'_{m+1\,I}|_{a_{m+1}+\sigma_m} \ge |h'_{m\,I}|_{a_{m+1}+\sigma_m} - |h'_{m+1\,I} - h'_{m\,I}|_{a_{m+1}+\sigma_m}\\
&&\ge \frac{1}{(2 - 2^{-m})c'_h} -  \frac{8 c_h c'_h s_m}{\sigma_m^2} \ge\frac{1}{(2 - 2^{-m})c'_h} - 288 \frac{2^{4n+10} c_h c'_h s_0}{a'^2_0} e^{-(c_s -4n-10)m - 2^m}\\
&&\ge \frac{1}{(2 - 2^{-m})c'_h} - \frac{2^m}{(2^{m+1}-1)(2^{m+2}-1)c'_h} = \frac{1}{(2 - 2^{-m-1})c'_h}.\\
\end{eqnarray*}
Consider the first inequality (\protect\ref{h_y}). For $c_s > c_{s\,3} = c_{s\,3}(c_h,c'_h,a_0,s_0)$
\begin{eqnarray*}
&&|h''_{m+1\,II}|_{a_{m+1}} \ge |h''_{m\,II}|_{a_{m+1}+\sigma_m} - |h''_{m+1\,II} - h''_{m\,II}|_{a_{m+1}}\\
&&\ge \frac{1}{(2 - 2^{-m})c'_h} -  \frac{8 c_h c'_h s_m}{\sigma_m^3} \ge\frac{1}{(2 - 2^{-m})c'_h} - 1728 \frac{2^{6n+15} c_h c'_h s_0}{a'^2_0} e^{-(c_s -6n-15)m - 2^m}\\
&&\ge \frac{1}{(2 - 2^{-m})c'_h} - \frac{2^m}{(2^{m+1}-1)(2^{m+2}-1)c'_h} = \frac{1}{(2 - 2^{-m-1})c'_h}.
\end{eqnarray*}
For the last two inequalities (\protect\ref{h_y}) let $c_s > c_{s\,4} = c_{s\,4}(c_h,c'_h,c''_h,a_0)$. Then
\begin{eqnarray*}
       |\sqrt \eps h''_{m+1\,Iy} - \sqrt \eps  h''_{m\,Iy}|_{a_{m+1}}
 &\le& \frac{8 n c_h c'_h s_m}
            {\sigma_m^3}\le 1728\frac{2^{6n+15} n c_h c'_h s_0}
                       {a'^2_0}e^{-(c_s - 6n-15)m - 2^m} \le c''_h 2^{-m-1} , \\
       |\eps h''_{m+1\,yy} - \eps h''_{m\,yy}|_{a_{m+1}}
 &\le& \eps \frac{8 n c_h c'_h s_m}
                       {\eps\sigma_m^3} \le c''_h 2^{-m-1}.
\end{eqnarray*}
For sufficently large
$c_s > \max\{c_{s\,1},c_{s\,2},c_{s\,3},c_{s\,4}\}$ the exponents
$$
  e^{-(c_s -2n-5)m - 2^m}, \quad
  e^{-(c_s -4n-10)m - 2^m}, \quad
  e^{-(c_s -6n-15)m - 2^m}
$$
are small and all inequalities (\ref{h_I}),(\protect\ref{h_y}) hold.

\subsection{Inequalities (\protect\ref{ang}),(\protect\ref{act})}

Note, that  for  $c_{\lambda} > c_{\lambda\,1}(c_N,c_s,b_0,a_0,s_0)$:
\begin{eqnarray*}
\frac{L_m s_m}{\delta_m} &\le& 3 \frac{2^{8n+11} c_N^{n+1} 2^{(4n+6)m}s_0 e^{-c_s m-2^m}}{c_{\lambda} b_0} \le \frac{1}{6} a'_0 2^{-(2n+5)(m+1)} = \sigma_m,\\
\frac{L_m s_m}{\sigma_m} &\le& 6\frac{2^{10n + 15} c_N^{n+1} 2^{(6n+10)m}s_0 e^{-c_s m-2^m}}{c_{\lambda} a'_0} \le \frac{1}{3} b_0 2^{-(m+1)} = \delta_m.\\
\end{eqnarray*}
This implies (\ref{ang}) and (\ref{act}).
\subsection{Inequalities (\protect\ref{cs}),
  (\protect\ref{AddStepang}),(\protect\ref{AddStepact})}

The first inequality in (\ref{cs}) holds for small $\eps$. Note, that for $c_s > c_{s\,5}(s_0,c'_h,a_0,b_0)$ we obtain
$$
s_m = s_0 e^{-c_s m - 2^m}  \le \frac{\sigma_m \delta_m}{4c'_h} = \frac{a'_0 b_0}{72c'_h}2^{-(2n+6)(m+1)}, \quad m \ge 1.
$$
By Proposition \ref{prop:Hn+1}:
$$
     |\widetilde S|_{a_m-4\sigma_m,b_m-2\delta_m}
 \le 16 \pi c'_h s'_m, \quad
     |\widetilde S'_{\hat\ph}|_{a_m-4\sigma_m,b_m-2\delta_m}
 \le 8c'_h s'_m.
$$
For $m=0$ inequalities (\ref{cs}), (\ref{AddStepang}), (\ref{AddStepact}) hold if $\eps$ is sufficently small. Consider $m \ge 1$. For  $c_s > c_{s\,6}(s_0,c'_h,a_0,b_0)$
\begin{eqnarray*}
       |\widetilde S'_{\hat I}|_{a_m-5\sigma_m,b_m-2\delta_m}
 &\le& 16 \pi c'_h \frac{s_m}{\sigma_m}
  \le 96 \pi c'_h \frac{s_0 2^{2n+5}}{a'_0}e^{-(c_s-2n-5)m - 2^m}
  \le \frac{b_0}{3}2^{-(m+1)} = \delta_m, \\
       |\sqrt \eps \widetilde S'_{\bar y}|_{a_m-5\sigma_m,
                                            b_m - 2\delta_m}
 &\le& \sqrt \eps 16 \pi c'_h \frac{s_m}{\sqrt \eps  \sigma_m}
  \le \delta_m,\\
  |\widetilde S'_{\hat\ph}|_{a_m-4\sigma_m,b_m-2\delta_m}
 &\le& 8c'_h  s_m \le 8c'_h  s_0 e^{-c_sm - 2^m} \le \frac{a'_0}{6}2^{-(2n+5)(m+1)} = \sigma_m.
\end{eqnarray*}

We choose the constants in the following way. Fix $c_s > \max(c_{s\,1},...c_{s\,6})$, $c_N > c_{N}(c_s)$ and $c_{\lambda} = c_{\lambda}(c_N,c_s) > \max(c_{\lambda\,0},c_{\lambda\,1})$, we obtain for $m \ge 0$ inequalities (\protect\ref{||m<}), (\protect\ref{ang}), (\protect\ref{act}) and for $m \ge 1$ we obtain (\protect\ref{h_I}), (\protect\ref{h_y}), (\protect\ref{cs}),(\protect\ref{AddStepang}),(\protect\ref{AddStepact}). Finally, for sufficently small $\eps$  inequalities (\protect\ref{h_I}), (\protect\ref{h_y}), (\protect\ref{cs}),(\protect\ref{AddStepang}),(\protect\ref{AddStepact}) hold for $m = 0$.

\section{Proofs}

\subsection{Proof of Proposition \protect\ref{prop:B}}

\begin{prop}
\label{prop:subset}
For any $K\in\mZ^{n+1}$, $|K|\le N_m$
$$
  Q_{K,m+1} \subset Q_{K,m}.
$$
\end{prop}

\begin{cor}
\label{cor:D}
$U_{a_{m+1}}(Q_{K,m+1}) \subset U_{a_{m}}(Q_{K,m})$ and equation (\ref{D=D}) which defines the domains $D_m$ can be represented as
$$
    D_{m+1}
  = D_m \setminus \cup_{N_{m-1}<|K|\le N_m}  U_{a_m}(Q_{K,m}) .
$$
\end{cor}
Consider the scaled frequency map
$$\tilde \nu_m: (y,I) \mapsto \Big(\Lambda'_y(y) + \eps h'_{m\,y}(y,I),h'_{m\,I}(y,I)\Big).$$
In comparison with (\ref{psi}) we remove the multiplier $\sqrt \eps$ at $h'_{m\,I}$.
It's Jacobi matrix equals
\begin{equation}
\label{Jacob}
  J_m(y,I) =
\begin{pmatrix}
\Lambda''_{yy} + \eps h''_{m\,yy} && h''_{m\,Iy}\\
\eps h''_{m\,yI} && h''_{m\,II}
\end{pmatrix}.
\end{equation}
\begin{prop}
\label{prop:Jacobt}
For some positive constants $\underline{C}_J$ and $\overline{C}_J$
\begin{equation}
\label{est:Jacob}
|{\rm det} J_m(y,I)|_{a_m} \le \overline{C}_J, \quad  |{\rm det} J_m(y,I)|^{-1}_{a_m} \le \underline{C}_J.
\end{equation}
\end{prop}
Estimates (\ref{est:Jacob}) imply the following inequality for measure of the domain $\tilde \nu_m(D_m)$:
$$
\mu \Big( \tilde \nu_m(D_m) \Big)  < \overline{C}_J \, \mu(D_m).
$$
Consider the vector $(\omega_y, \sqrt \eps \omega_I),\quad \omega_y \in {\mR}^n,\quad \omega_I \in \mR$ and set
\begin{eqnarray}
\label{Qomega}
&     Q^{\omega}_{K,m}
   =  \big\{ (\omega_y, \omega_I) \in \tilde \nu_m(D_m) :
             \big| \langle \omega_y ,k\rangle + \sqrt \eps \omega_I k_0
             \big|
         \le \lambda_{\bj(|K|)} (1+2^{-m-1}) \sqrt \eps
      \big\}.
\end{eqnarray}
The set  $Q^{\omega}_{K,m} \subset \mR^{n+1}$ is a strip between two planes
$$ \frac{\langle \omega_y ,k\rangle}{\sqrt{\langle k ,k\rangle + \eps k_0^2}} + \frac{\sqrt \eps \omega_I k_0}{\sqrt{\langle k ,k\rangle + \eps k_0^2}}= \pm \frac{\sqrt \eps \lambda_{\bj(|K|)} (1+2^{-m-1}) }{\sqrt{\langle k ,k\rangle + \eps k_0^2}}.$$
Using (\ref{cutoff}) we have, that $\langle k,k\rangle \ge 1$ and the the distance between the plains is not more than $4 \lambda_{\bj(|K|)} \sqrt \eps$.
The measure estimates are
\begin{eqnarray*}
\mu\Big(Q^{\omega}_{K,m} \Big) &\le& 4\lambda_{\bj(|K|)}\sqrt \eps\, \overline{C}_J \,C_{D},\\
 \mu\Big( Q_{K,m} \cap D_m \Big) &\le& \mu\Big( \tilde \nu_m^{-1}(Q^{\omega}_{K,m}) \Big) \le 4\lambda_{\bj(|K|)}\sqrt \eps\, \underline{C}_J \, \overline{C}_J \,C_{D},
\end{eqnarray*}
where $C_{D}$ depends on diameter and dimension of $D_0$.

Consider estimates for the measure of $D_m \cap U_{a_m}(Q_{K,m})$. Let $|y'| \le \sqrt \eps \sigma_m$, $|I'| \le	\sigma_m$. Than
\begin{eqnarray*}
&&|\langle \nu_{m}(y + y',I + I') - \nu_{m}(y,I),K \rangle| \le n |\Lambda''_{yy}| |K|\sqrt \eps \sigma_m + n |\eps h''_{yy}| |K|\sqrt \eps \sigma_m\\
&& + n |\sqrt \eps h''_{yI}| |K|\sqrt \eps \sigma_m +  n |\sqrt \eps h''_{yI}| |K|\sqrt \eps \sigma_m + n |\sqrt \eps h''_{yI}| |K|\sqrt \eps \sigma_m\\
&&\le \sqrt \eps C_{\psi} N_{m} \sigma_m,
\end{eqnarray*}
where $C_{\psi} = C_{\psi}(c_{\Lambda},c'_h,c''_h,n)$.

Consider the extension of  $Q^{\omega}_{K,m}$
\begin{eqnarray}
\label{Qomega+}
&     Q^{\omega+}_{K,m}
   =  \big\{ (\omega_y, \omega_I) \in \tilde \nu_{m}(D_m) :
             \big| \langle \omega_y ,k\rangle + \sqrt \eps \omega_I k_0
             \big|
         \le (2 \lambda_{\bj(|K|)} + C_{\psi} N_m \sigma_m) \sqrt \eps
      \big\}.
\end{eqnarray}
Note that $\Big(D_m \cap U_{a_m}(Q_{K,m})\Big) \subset \tilde\nu^{-1}_m(Q^{\omega+}_{K,m})$. Finally
\begin{equation}
\label{muest}
\mu\Big(D_m \cap U_{a_m}(Q_{K,m}\Big) \le \mu(\tilde \nu^{-1}_m(Q^{\omega+}_{K,m})) \le 4(\lambda_{\bj(|K|)} + C_{\psi} N_m \sigma_m) \sqrt \eps\, \underline{C}_J \, \overline{C}_J C_{D}.
\end{equation}
We have $( D_m \setminus D_{m+1}) \subset \cup_{N_{m-1}<|K|\le N_m}\Big(D_m \cap U_{a_m}(Q_{K,m})\Big)$. Using (\ref{muest}) we obtain
\begin{eqnarray*}
\mu \Big( D_0 \setminus D_{m+1} \Big) &\le& \sum_{i = 0}^{m} \mu \Big( \cup_{N_{i-1}<|K|\le N_i} \Big( D_i \cap U_{a_i}(Q_{K,i})\Big) \Big) \le \\
		&\le& 4 \sqrt \eps \, \underline{C}_J \, \overline{C}_J C_{D}  \sum_{i = 0}^{m} N_i^{n+1}(\lambda_{i} + C_{\psi} N_i \sigma_i).
\end{eqnarray*}

The proposition holds for $c_{\mu} = 4\underline{C}_J \, \overline{C}_J C_{D}  \sum_{i = 0}^{m} N_i^{n+1}(\lambda_{i} + C_{\psi} N_i \sigma_i)$. To finish the proof we need to check
\begin{equation}
\label{musum}
\sum_{i = 0}^{+\infty} N_i^{n+1}(\lambda_{i} + C_{\psi} N_i \sigma_i) < +\infty.
\end{equation}
Using (\ref{defsigmadelta}) and (\ref{defsNLambda}) we obtain
\begin{eqnarray*}
\sum_{i = 0}^{+\infty} N_{i}^{n+1}(\lambda_{i} + C_{\psi} N_{i} \sigma_i) &=& \sum_{i = 0}^{+\infty} (c_N^{n+1} c_{\lambda} 2^{-i} + \frac{1}{6} a_0 C_{\psi} c_N^{n+2} 2^{-i}) =\\
&=&  2 c_N^{n+1} c_{\lambda} + \frac{1}{3} a_0 C_{\psi} c_N^{n+2} < +\infty.
\end{eqnarray*}
\qed

\subsection{Proof of Proposition \protect\ref{prop:homolo}}

Solution of equation (\ref{homologic}) has the form
\begin{equation}
\label{formal}
    S
  = \sum_{|K|\le N_{m}, \, k\ne 0} \!
      S^K\, e^{i\langle k,x \rangle + ik_0\ph}, \quad
    S^K
  =  \frac{i \sqrt \eps \, V_m^K(\hat y,\hat I)}
          {\langle\nu_m(\hat y,\hat I), K \rangle } ,
\end{equation}
where $\nu_m(\hat y,\hat I)$ is determined by (\ref{psi}).

The first inequality (\ref{||m<}) means that
$$
      |V_m^K e^{i\langle k,x\rangle + ik_0\ph}|_{a_m,b_m-\delta}
  \le 2s_m e^{-(|k_1| + |k_2| + ... + |k_0|) \delta}, \qquad
      0\le\delta\le b_m.
$$
Then
$$
      |S|_{a_m,b_m-\delta_m}
  \le \sum_{j=0}^{m} \sum_{N_{j-1} < |K| \le N_j,\, k\ne 0}
      \frac{2s_m}{\lambda_j}
      \, e^{-\delta_m (n+1) N_{j-1}}
  \le L_m s_m .
$$\qed

\subsection{Proof of Proposition \protect\ref{prop:RRRR}}

In this section for brevity we write $V,h$ instead of $V_m,h_m$ and
$a,b,\sigma,\delta,N,L,s,\tilde s$ instead of
$a_m$, $b_m$, $\sigma_m$, $\delta_m$, $N_{m}$, $L_m$,
$s_m$, $\tilde s_m$.

The function $\widetilde V_m = \tilde v_m + \tilde u_m$ can be presented in the form
\begin{eqnarray}
\label{V=RRRRR}
 &&  \widetilde V_m(\hat y,\hat I,\hat x,\hat\ph)
  = R_1 + R_2 + R_3 + R_4 + R_5 ,  \\
\nonumber
     R_1
 &=& \frac{1}{\eps}
     \Big( \Lambda(y) - \Lambda(\hat y)
           - \langle \Lambda^\prime_y(\hat y),
             \sqrt\eps  S^\prime_x \rangle
     \Big), \\
\nonumber
     R_2
 &=& h(y,I) - h(\hat y,\hat I)
      - \big\langle h'_y(\hat y,\hat I),
             \sqrt\eps S^\prime_x \big\rangle
      - h'_I(\hat y,\hat I) \,
             S^\prime_\ph , \\
\nonumber
     R_3
 &=& V(y,I,x,\ph) - V(\hat y,\hat I,x,\ph) , \\
\nonumber
     R_4
 &=& V(\hat y,\hat I,x,\ph)
      - \langle V \rangle_x (\hat y,\hat I,\ph)
      - \Pi_{N} V(\hat y,\hat I,x,\ph) , \\
\nonumber
     R_5
 &=& \langle V \rangle_x (\hat y,\hat I,\ph)
    - \langle V \rangle_x (\hat y,\hat I,\hat\ph) .
\end{eqnarray}

By Proposition (\ref{prop:homolo}) the first term in (\ref{V=RRRRR}) satisfies the estimate
$$
       |R_1|_{a - \sigma, b - 2\delta}
  \le  \frac12 |\Lambda''_{yy}|_a \,
       |S^\prime_x|^2_{a - \sigma, b - 2\delta}
  \le \frac{c_\Lambda }{2} \,
              \Big(\frac{L s}{\delta}\Big)^2 .
$$
To estimate the second one we use (\ref{h_I}),(\ref{h_y}):
\begin{eqnarray*}
       |R_2|_{a - \sigma, b - 2\delta}
 &\le& \frac12
     \Big( |\eps h''_{yy}|_a \,
        |S'_x|^2_{a - \sigma, b - 2\delta} + 2n |\sqrt\eps h''_{yI}|_a \,
    |S'_x S'_\ph|_{a - \sigma, b - 2\delta} + |h''_{II}|_a \,
             |S'_\ph|^2_{a - \sigma, b - 2\delta}
       \Big) \le \\
 &\le&
           c''_h  n \,
           \Big(\frac{L s}{\delta}\Big)^2
          + 2 n c''_h \, \Big(\frac{ L s}{\delta}\Big)^2
          + c'_h \,
            \Big(\frac{L s}{\delta}\Big)^2 =  (c'_h + 3n c''_h)\Big(\frac{L s}{\delta}\Big)^2
\end{eqnarray*}
The third term is estimated by (\ref{||m<}):
\begin{eqnarray*}
       |R_3|_{a - \sigma, b - 2\delta}
 &\le&  \sqrt \eps n |V'_y|_{a - \sigma,b - \delta}
         |S'_x|_{a - \sigma,b - 2\delta}
      + |V'_I|_{a - \sigma,b - \delta}
         |S'_\ph|_{a - \sigma,b - 2\delta} \le \\
 &\le& (n + 1) \,
        \frac{L s^2}{\sigma \delta}.
\end{eqnarray*}
By (\ref{PIN})
$$
      |R_4|_{a - \sigma, b - 2\delta}
  \le \frac{C}{\delta} \Big( N + \frac{1}{\delta} \Big)^n
           e^{- N \delta} \, s .
$$
Finally
$$
      |R_5|_{a - \sigma, b - 2\delta}
    \le |V'_{\ph}|_{a - \sigma, b - \delta}
         |S'_I|_{a - \sigma,b - \delta}
    \le  \frac{L s^2}{\sigma\delta} .
$$
Note, that $3nc''_h \le \frac{c_{\Lambda}}{2}$. Therefore
$$
     |\tilde u_m + \tilde v_m |_{a - \sigma,b - 2\delta}
 \le \tilde s,
$$
where
$$
       \tilde s
   =   (c_{\Lambda} + c'_h)\Big(\frac{L s}{\delta}\Big)^2
     + (n+2) \frac{L s^2}{\sigma\delta} + \frac{C}{\delta} \Big(N + \frac{1}{\delta}\Big)^n
       e^{-N \delta}\, s .
$$
\qed

\subsection{Proof of Proposition \ref{prop:Hn+1}}

By (\ref{0}), (\ref{Lam-1}), (\ref{||m<}), (\ref{h_I}), and (\ref{|v|}), for any $m\ge 1$ we have:
\begin{eqnarray}
\label{hhv}
&\!\!\!\!\!\!\!\!\!\!\!\!
  |h_0|_{a_0} \le c_h, \quad
  |h'_{0I}|^{-1}_{a_0} \le c'_h, \quad
  |v_0|_{a_0,b_0} \le \sqrt\eps\bc, \quad
  |h_m|_{a_m} \le 2c_h, \quad
  |h'_{mI}|_{a_m}^{-1} \le 2c'_h, & \\
\label{hmhmvm}
&   |v_m|_{a_m-\sigma_m,b_m-2\delta_m}
  \le s_m, \quad   |\tilde v_m + \tilde u_m|_{a_m-\sigma_m,b_m-2\delta_m}
  \le \tilde s_m.&
\end{eqnarray}
Applying Lemma \ref{lem:aa} to equation (\ref{hm+1}), we obtain:
\begin{eqnarray}
\label{|tildeS|}
&\!\!\!\!\!\!\!\!\!\!\!\! \!\!\!\!\!\!\!\!\!\!\!\! \!\!\!\!\!\!\!\!\!\!\!\!   |\widetilde S'_\ph|_{a_m - 4\sigma_m,b_m - 2\delta_m}
  \le 8c'_h |v_m|_{a_m - \sigma_m,b_m - 2\delta_m}, \quad
      |\widetilde S|_{a_m - 4\sigma_m,b_m - 2\delta_m}
  \le 16 \pi c'_h |v_m|_{a_m - \sigma_m,b_m - 2\delta_m}, & \\
\label{h-h}
&\displaystyle
      \widetilde S'_x = 0, \quad
      |h_m - h_{m+1}|_{a_m - 4\sigma_m,b_m - 2\delta_m}
  \le 8c_h c'_h
       \frac{|v_m|_{a_m - \sigma_m,b_m - 2\delta_m}}
            {\sigma_m} .   &
\end{eqnarray}
Then estimates (\ref{|h-h|}) follow from (\ref{hhv}), (\ref{hmhmvm}), and (\ref{h-h}).

We have:
$$
      v_{m+1}(\bar y,\bar I,\bar\ph)
    + u_{m+1}(\bar y,\bar I,\bar x,\bar\ph)
  =  \tilde v_{m}(\bar y,\bar I + \widetilde S'_{\hat\ph},\bar\ph - \widetilde S'_{\hat I}) +  \tilde u_m \big(\bar y,\bar I + \widetilde S'_{\hat\ph},
                       \bar x - \sqrt\eps\widetilde S'_{\bar y},
                        \bar\ph - \widetilde S'_{\hat I}
                 \big) .
$$
Since by (\ref{|tildeS|})
\begin{eqnarray*}
     |\widetilde S'_{\hat\ph}|_{a_{m+1},b_{m+1}+2\delta_m}
&\le& 8 c'_h s_m
 \le  \sigma_m, \\
      |\sqrt\eps \widetilde S'_{\bar y}|_{a_{m+1},b_{m+1}+\delta_m}
&\le& \sqrt \eps \frac{16 \pi c'_h  s_m}{\sigma_m} \le\delta_m, \\
      |\widetilde S'_{\bar I}|_{a_{m+1},b_{m+1}+\delta_m}
&\le& \frac{16 \pi c'_h s_m}{\sigma_m} \le\delta_m,
\end{eqnarray*}
we have:
$$
      |v_{m+1} + u_{m+1}|_{a_{m+1},b_{m+1}}
  \le \tilde s_m \le s_{m+1}.
$$
\qed

\subsection{Proof of Proposition \ref{prop:L}}

By equations (\ref{L}) and (\ref{defsNLambda}) we have

\begin{eqnarray*}
 L_m \le \sum_{j=0}^{m}\frac{2 N_j^{n+1}}{\lambda_j} = \frac{2 c_N^{n+1}}{c_{\lambda}}\sum_{j=0}^{m} 2^{(4n+5)j} \le \frac{c_N^{n+1}}{c_{\lambda}}2^{(4n+5)(m+2)}.
\end{eqnarray*}
\qed

\subsection{Proof of Proposition \protect\ref{prop:subset}}
A point $(y,I)\in Q_{K,m+1}$ if$$
      \big|\big\langle \nu_{m+1}(y,I), K
      \big\rangle\big|
  \le \lambda_{\bj(|K|)} (1 + 2^{-m-2}) \sqrt \eps.
$$
Our aim is to show that $(y,I)\in Q_{K,m}$, i.e.
\begin{equation}
\label{<nu,K>}
      \big|\big\langle \nu_m(y,I), K \big\rangle\big|
  \le \lambda_{\bj(|K|)} (1 + 2^{-m-1}) \sqrt \eps.
\end{equation}
We have the inequality
\begin{eqnarray*}
& \big| \langle \nu_m(y,I), K \rangle
       - \langle \nu_{m+1}(y,I), K \rangle \big|
 \le (n+1)W\, |K| , & \\
   & W
  = |\nu_m(y,I) - \nu_{m+1}(y,I)| .
\end{eqnarray*}
By (\ref{|h-h|})
$$
       W
  \le  \eps |h'_{m+1\, y} - h'_{m\, y}|_{a_{m+1},b_{m+1}}
      + \sqrt\eps |h'_{m+1\, I} - h'_{m\, I}|_{a_{m+1},b_{m+1}}
  \le  \sqrt \eps \frac{16 c_h c'_h s'_m}{\sigma^2_m}.
$$
Therefore for any $|K|\le N_m$ we have the estimate
\begin{eqnarray*}
     (n+1)W |K|
 \le (n+1) \frac{16 c_h c'_h s'_m}{\sigma^2_m}N_m,
\end{eqnarray*}
where $s'_0 = \sqrt\eps \bc$ and $s'_m = s_m$ for all $m\ge 1$.

It remains to check the estimate
\begin{equation*}
 (n+1) \frac{16 c_h c'_h s'_m}{\sigma^2_m}N_m  \le \lambda_{m} 2^{-m-2}.
\end{equation*}
Let $c_{\lambda} > c_{\lambda}(c_s,c_N)$. For $m=0$ we have
\begin{eqnarray*}
(n+1)\frac{16 c_h c'_h s'_0}{\sigma^2_0}N_0 \le \lambda_{0} 2^{-2}.
\end{eqnarray*}
and for $m \ge 1$
\begin{eqnarray*}
(n+1) \frac{36 c_h c'_h s_0 c_N}{a'^2_{0}}2^{4n+14}2^{(4n + 12)m}e^{-c_sm - 2^m}  \le c_{\lambda} 2^{-(2n+3)m}.
\end{eqnarray*}
\qed

\subsection{Proof of Proposition \protect\ref{est:Jacob}}
Suppose that the arguments of functions $\Lambda$, $h$ lie in $U_{a_m}(D_m)$. Let us expand the Jacobian $J_m$ with respect to the last column
\begin{eqnarray*}
\det J_m(y,I) &=& {\rm det}  \begin{pmatrix}
\Lambda''_{yy} + \eps h''_{m\,yy} && \sqrt \eps  h''_{m\,Iy}\\
\sqrt \eps h''_{m\,yI} && h''_{m\, II}
\end{pmatrix} =\\
&=& {\rm det} \Big( \Lambda''_{yy} + \eps h''_{m\,yy} \Big) h''_{m\,II} + \sum_{i=1}^{n} (-1)^{n+1+i} \sqrt \eps  h''_{m\,Iyi} M_{i,n+1}.
\end{eqnarray*}
Here $M_{i,n+1}$ is the $(i,n+1)$ minor matrix of $J_m(y,I)$. Using (\ref{Lam-1}), (\ref{h0''}), (\ref{h_I}) and (\ref{h_y}) we obtain
$$
|\det J_m(y,I)| \le  n! \Big(|\Lambda''_{yy}| + |\eps h''_{m\,yy}|\Big)^n |h''_{m\,II}| + n n! |\sqrt \eps  h''_{m\,Iy}|^2 \Big(|\Lambda''_{yy}| + |\eps h''_{m\,yy}|\Big)^{n-1} \le \overline{C}_J,
$$
where $\underline{C}_J$ is some constant depending on $c_{\lambda},c'_h,c''_h$ and $n$.

Return to estimate for the Jacobian. For sufficiently small $c''_h$
$$c''_h \le c''_{h\,1} = \min \Big( \frac{\underline{c}_{\Lambda}}{2^{n+1}nn!c_{\Lambda}^{n-1}}, \frac{c_{\Lambda}}{2}\Big), \quad |\eps h''_{m\,yy}| \le 2 c''_h$$
 we have
$$
\Big|\det \Big( \Lambda''_{yy} + \eps h''_{m\,yy}\Big)\Big| > |\det \Lambda''_{yy}| - n n! |\eps h''_{m\,yy}| \Big( |\Lambda''_{yy}| + |\eps h''_{m\,yy}| \Big)^{n-1} \ge  \frac{1}{2}\underline{c}_{\Lambda},
$$
Note, that $|h''_{m\,II}| \ge \frac{1}{2c'_h}$ and for
$$c''_h \le c''_{h\,2} = \min \Big( \frac{\underline{c}_{\Lambda}}{2^{n+4}nn!c^{n-1}_{\Lambda}c'_h}, c''_{h\,1}\Big)$$
we have
$$\Big|\sum_{i=1}^{n} (-1)^{n+1+i} \sqrt \eps  h''_{m\,Iyi} M_{i,n+1}\Big| \le \frac{\underline{c}_{\Lambda}}{8 c'_h}.$$
Finally
$$
|\det J(y,I)|_{a_m} \ge \frac{\underline{c}_{\Lambda}}{8 c'_h} = \underline{C}_{J}^{-1}.
$$
\qed

\section{Further technical statements}

\subsection{Lemma on a cut off}

\begin{lem}
\label{lem:cutoff}
For any real-analytic function $f$ on $U_b(\mT^{m+1})$ and any $\delta\in (0,b)$
$$
    \big| f - \langle f\rangle_x - \Pi_N f
    \big|_{b-\delta}
  \le \frac{C}{\delta} \Big(N + \frac{1}{\delta}\Big)^n
         e^{- N\delta} |f|_b .
$$
where the constant $C$ depends only on $n$.
\end{lem}

{\it Proof}. The Fourier coefficients (\ref{cutoff}) satisfy the inequalities
$$
  |f^K| \le e^{-b |K|} |f|_b .
$$
Then the equation
$$
    f - \langle f\rangle_x - \Pi_N f
  = \sum_{|K|>N,\, k\ne 0} f^K e^{i\langle k,x\rangle + k_0\ph}
$$
implies
$$
      |f - \langle f\rangle_x - \Pi_N f|_{b-\delta}
  \le |f|_b \sum_{|K|>N,\, k\ne 0} e^{-\delta |K|} .
$$
The sum in the right-hand side does not exceed
$$
      c_1 \int_{x\in\mR^{n+1},\, |x|>N} e^{-\delta |x|} \, dx
  \le \frac{c_2}{\delta^{n+1}}
        \int_{\delta N}^\infty s^n e^{-s} \, ds
  \le \frac{c_3}{\delta^{n+1}} (1 + \delta N)^n e^{-\delta N},
$$
where $c_1,c_2,c_3$ depend only on $n$.

\subsection{Lemma on the action-angle variables}

\begin{lem}
\label{lem:aa}
Let $h$ and $v$ be real-analytic functions, defined in complex
neighborhoods of $[-\alpha,\alpha]$ and $[-\alpha,\alpha]\times\mT$
respectively. Let the canonical change
$(I,\ph\bmod 2\pi)\mapsto (\bar I,\bar\ph\bmod 2\pi)$, determined by the generating function $\bar I\ph + S(\bar I,\ph)$,
$\langle S\rangle_\ph = 0$, be such that
\begin{eqnarray}
\label{h+v}
 &  h(I) + v(I,\ph) = h_*(\bar I), & \\
\label{h'>}
 & |h|_a \le c, \quad
   |h'|_a^{-1} \le c', \quad
   |v|_{a,b} \le \frac{\sigma}{2c'}, \quad
   0 < \sigma < a/2 &
\end{eqnarray}
Then
\begin{equation}
\label{aa}
      |S'_\ph|_{a-3\sigma,b}
  \le 4 c' |v|_{a,b}, \quad
      |S|_{a-3\sigma,b}
  \le 8 \pi c' |v|_{a,b}, \quad
      |h - h_*|_{a-3\sigma}
  \le 2 c c' \frac{|v|_{a,b}}{\sigma}.
\end{equation}
\end{lem}

{\it Proof}. Equation that determines $\bar I$ is well-known:
\begin{equation}
\label{barI}
    \bar I(r)
  = \frac 1{2\pi} \int_0^{2\pi} I(r,\ph)\, d\ph.
\end{equation}
Here $I(r,\ph)$ is the solution of the equation
$$
  h(I) + v(I,\ph) = h(r).
$$
We use $r$ as a constant which fixes the energy $h(r)$.

By Lemma \ref{lem:hvh} the function $I = I(r,\ph)$ is as follows:
\begin{equation}
\label{I=r+f}
  I = r + f(r,\ph), \quad
  |f|_{a-2\sigma,b} \le 2c'|v|_{a,b}.
\end{equation}
Moreover,
\begin{equation}
\label{implies}
  I\in U_{a-2\sigma-\delta}([-\alpha,\alpha])\quad
  \mbox{implies}\quad
  r\in U_{a-\delta}([-\alpha,\alpha])\quad
  \mbox{for any } \delta\in [2\sigma,a-\sigma].
\end{equation}

Equations (\ref{barI}) and (\ref{I=r+f}) imply
\begin{eqnarray}
\label{I-barI}
&  I(r,\ph) - \bar I(r)
 = f(r,\ph) - \langle f \rangle_\ph(r), & \\
\label{rI}
&  r = r(\bar I), \quad
   I - \bar I = S'_\ph(\bar I,\ph). &
\end{eqnarray}
Combining (\ref{I-barI}) and (\ref{rI}), we obtain:
$$
    S(\bar I,\ph)
  = \int_0^{\ph} \big( I(r,\ph) - \bar I \big)\, d\ph
  = \int_0^{\ph}
                  \big( f(r,\ph) - \langle f \rangle_\ph(r) \big)
                    \, d\ph.
$$
Therefore
$$
   \quad |S'_\ph|_{a-3\sigma,b} \le 4 c' |v|_{a,b},\quad |S|_{a-3\sigma,b} \le 8 \pi c' |v|_{a,b}.
$$

By using the equation $h(r)=h_*(\hat I)$, we have:
\begin{eqnarray*}
       |h(I) - h_*(I)|_{a-3\sigma}
 &\le& |h(r + f(r,\ph)) - h(r)|_{a-2\sigma}
  \le  |h'|_{a-\sigma} |f|_{a-2\sigma,b}
  \le  \frac{c}{\sigma}\, 2c' |v|_{a,b}.
\end{eqnarray*}
\qed

\subsection{A version of the implicit function theorem}

\begin{lem}
\label{lem:hvh}
 Let the real-analytic functions $h,v$, defined in a complex
neighborhood of the interval $\II\subset\mR$, satisfy the estimates
\begin{equation}
\label{|h|,|v|}
    |h'|_a^{-1} \le c', \quad
    |v|_a \le \frac{\sigma}{2c'}, \qquad
    0 < \sigma < \frac a2.
\end{equation}
Then the equation
\begin{equation}
\label{hvh}
   h(I) + v(I) = h(r), \qquad
   I\in U_{a-\sigma}(\II)
\end{equation}
implies
\begin{equation}
\label{|f|}
  I = r + f(r), \qquad
  |f|_{a-2\sigma} \le 2c' |v|_a \le \sigma,
\end{equation}
where $f(r)$ is the real-analytic function, $r \in U_{a-2\sigma}(\II)$.
\end{lem}

{\it Proof}. Applying the map $h^{-1}$ to (\ref{hvh}), we get:
$$
  I + u(I) = r, \qquad
  u(I) = h^{-1} \big( h(I) + v(I) \big) - I.
$$
If $|v|_a \le \sigma / c'$, the function $u$ is defined in a
complex neighborhood of $\II$ and admits the estimate
$$
      |u|_{a-\sigma}
  \le |h'|_a^{-1} |v|_a
  \le c' |v|_a .
$$

The function $I = I(r)$, defined by (\ref{|f|}), is a fixed point of the operator
$$
  I(r) \mapsto \Phi(I(r),r) = r - u(I(r)).
$$
This operator is contracting with respect to the norm $|\cdot|_{a-2\sigma}$ because by (\ref{|h|,|v|})
$$
     |\Phi'_I|_{a-2\sigma}
  =  |u'_I|_{a-2\sigma}
 \le \frac{c' |v|_a}{\sigma}
 \le \frac12.
$$
Therefore $I-r = f(r)$, where
$|f(r)|_{a-2\sigma} \le 2 |u|_{a-\sigma} \le 2c' |v|_a$.
\qed

\section*{Acknowledgements}

The work was supported by grants RFBR 13-01-00251, 05-01-01119, 13-01-12462, State contract no. 8223 ``Dynamic instability and catastrophe'' (RF) and  RF Program for the State Support of Leading Scientific Schools  NSh-2519.2012.1.

\end{document}